\documentclass[12pt]{article}
\usepackage{amsfonts}
\usepackage{amssymb}
\usepackage{amsmath}
\usepackage{color}
=0pt
\def\sqr#1#2{{\vcenter{\vbox{\hrule height.#2pt
              \hbox{\vrule width.#2pt height#1pt \kern#1pt \vrule
width.#2pt}
              \hrule height.#2pt}}}}
\def\signed #1{{\unskip\nobreak\hfil\penalty50
              \hskip2em\hbox{}\nobreak\hfil#1
              \parfillskip=0pt \finalhyphendemerits=0 \par}}
\def\endpf{\signed {$\sqr69$}}
\def\dbE{{\mathbb{E}}}
\def\dbF{{\mathbb{F}}}

\def\dbH{{\mathbb{H}}}

\def\dbN{{\mathbb{N}}}

\def\dbP{{\mathbb{P}}}

\def\dbR{{\mathbb{R}}}
\def\dbS{{\mathbb{S}}}

\def\Om{\Omega}
\def\om{\omega}
\def\a{\alpha}

\def\g{\gamma}
\def\d{\delta}
\def\e{\varepsilon}

\def\k{\kappa}
\def\l{\lambda}

\def\si{\sigma}

\def\f{\varphi}
\def\th{\theta}

\def\i{\infty}
\def\3n{\negthinspace \negthinspace \negthinspace }
\def\2n{\negthinspace \negthinspace }
\def\1n{\negthinspace }
\def\ns{\noalign{\smallskip} }

\def\ds{\displaystyle}
\def\G{\Gamma}
\def\D{\Delta}
\def\Th{\Theta}
\def\L{\Lambda}

\def\cA{{\cal A}}
\def\cB{{\cal B}}
\def\cC{{\cal C}}
\def\cD{{\cal D}}

\def\cF{{\cal F}}
\def\cG{{\cal G}}

\def\cJ{{\cal J}}

\def\cL{{\cal L}}

\def\cR{{\cal R}}
\def\cS{{\cal S}}

\def\cU{{\cal U}}

\def\cX{{\cal X}}

\def\cl{{\cal l}}
\def\mE{{\mathbb{E}}}

\def\no{\noindent}

\def\ms{\medskip}
\def\bs{\bigskip}
\def\q{\quad}
\def\qq{\qquad}
\def\hb{\hbox}

\def\Ra{\mathop{\Rightarrow}}

\def\lan{\mathop{\langle}}
\def\ran{\mathop{\rangle}}

\def\wt{\widetilde}
\def\cd{\cdot}
\def\cds{\cdots}

\def\ae{\hbox{\rm a.e.{ }}}

\def\span{\hbox{\rm span$\,$}}

\def\cl{\overline}

\def\deq{\mathop{\buildrel\D\over=}}

\def\({\Big (}
\def\){\Big )}
\def\[{\Big[}
\def\]{\Big]}

\def\={\buildrel \triangle \over =}

\def\lan{\big\langle}
\def\ran{\big\rangle}

\def\h{\hat}

\def\resp{{\it resp. }}

\def\be{\begin{equation}}
\def\bel{\begin{equation}\label}
\def\ee{\end{equation}}
\def\bea{\begin{eqnarray}}
\def\eea{\end{eqnarray}}
\def\bt{\begin{theorem}}
\def\et{\end{theorem}}
\def\bc{\begin{corollary}}
\def\ec{\end{corollary}}
\def\bl{\begin{lemma}}
\def\el{\end{lemma}}
\def\bp{\begin{proposition}}
\def\ep{\end{proposition}}
\def\br{\begin{remark}}
\def\er{\end{remark}}
\def\ba{\begin{array}}
\def\ea{\end{array}}
\def\bd{\begin{definition}}
\def\ed{\end{definition}}

\newtheorem{lemma}{Lemma}[section]
\newtheorem{remark}{Remark}[section]

\newtheorem{theorem}{Theorem}[section]
\newtheorem{corollary}{Corollary}[section]

\newtheorem{definition}{Definition}[section]
\newtheorem{proposition}{Proposition}[section]

\allowdisplaybreaks

\makeatletter
   
   \@addtoreset{equation}{section}
\makeatother

\begin{document}

\title{\bf Well-posedness of Stochastic Riccati
Equations and Closed-Loop Solvability for Stochastic Linear Quadratic Optimal Control
Problems}
\author{Qi L\"u\thanks{School of Mathematics, Sichuan
University, Chengdu, 610064, China. This work is
supported the NSF of China under grant 11471231,
the Fundamental Research Funds for the Central
Universities in China under grant 2015SCU04A02,
the NSFC-CNRS Joint Research Project under grant
11711530142 and Grant MTM2014-52347 of the
MICINN, Spain.}}

\date{}
\maketitle

\begin{abstract}
We study the closed-loop solvability of a
stochastic linear quadratic optimal control
problem for systems governed by stochastic
evolution equations. This solvability is
established by means of solvability of the
corresponding Riccati equation, which is implied
by the uniform convexity of the quadratic cost
functional. At last, conditions ensuring the
uniform convexity of the cost functional are
discussed.

\end{abstract}

\no{\bf 2010 Mathematics Subject
Classification}. 93E20, 49N10, 49N35.

\bs

\no{\bf Key Words}. stochastic linear quadratic
control problem, stochastic evolution equation,
closed-loop solvability,  Riccati equation.

\section{Introduction}

Let $(\Om,\cF,\dbF,\dbP)$ be a complete filtered
probability space on which a standard
one-dimensional Brownian motion
$W(\cd)=\{W(t)\}_{t\geq 0}$ is defined, and
$\dbF=\{\cF_t\}_{t\geq0}$ is the natural
filtration of $W(\cd)$ augmented by all the
$\dbP$-null sets in $\cF$.

Let $T>0$ be a fixed time horizon. For any
$t\in[0,T)$ and Banach space $\dbH$, let
$$\begin{array}{ll}\ds
L^2_{\cF_t}(\Om;\dbH)=\Big\{\xi:\Om\to\dbH\bigm|\xi\hb{
is
$\cF_t$-strongly measurable, }\dbE|\xi|_{\dbH}^2<\i\Big\},\\
\ns\ds
L_\dbF^2(t,T;\dbH)=\Big\{\f:[t,T]\times\Om\to\dbH\bigm|\f(\cd)\hb{
is
$\dbF$-progressively measurable},\\
\ns\ds\hspace{6.4cm}\dbE\int^T_t|\f(s)|_{\dbH}^2ds<\i\Big\},
\end{array}
$$
$$\begin{array}{ll}
\ds
C_\dbF([t,T];L^2(\Om;\dbH))=\Big\{\f:[t,T]\times\Om\to\dbH\bigm|\f(\cd)\hb{
is $\dbF$-adapted, }\f:[t,T]\to L^2_{\cF_T}(\Om;\dbH)\\
\ns\ds\hspace{7.694cm}\hb{is continuous,
 }\Big\},
\end{array}$$
$$\begin{array}{ll}
\ds
L^2_\dbF(\Om;L^1(t,T;\dbH))=\Big\{\f:[t,T]\times
\Om\to\dbH\bigm|\f(\cd)\hb{ is
$\dbF$-progressively measurable},\\
\ns\ds\hspace{7.64cm}
\dbE\(\int_t^T|\f(s)|_{\dbH}ds\)^2<\i\Big\}.
\end{array}$$

Let $\dbH_1$ and  $\dbH_2$ be two Banach spaces.
Denote by $\cL(\dbH_1;\dbH_2)$(\resp
$\cL(\dbH_1)$) the set of all bounded linear
operators from $\dbH_1$ to $\dbH_2$(\resp
$\dbH_1$). If $\dbH$ is a Hilbert space, then we
set
$$
\dbS(\dbH)\=\{F\in \cL(\dbH)|\, F \mbox{ is
symmetric }\},
$$
and
$$
\cl{\dbS_+}(\dbH)\=\big\{F\in \dbS(\dbH)\big|\,
\lan F\xi,\xi\ran\geq 0,\q  \forall \xi\in
H\big\}.
$$
Here and in what follows, for simplicity of
notations, when there is no confusion, \emph{we
shall use $\langle \cd\,,\cd\rangle$ for inner
products in possibly different Hilbert spaces.}

For any interval $[t_1,t_2]\subset [0,+\infty)$,
denote by $C([t_1,t_2];\dbS(\dbH))$ the set of
all continuous mappings from $[t_1,t_2]$ to
$\dbS(\dbH)$, which is a Banach space with the
norm
$$
|F|_{C([t_1,t_2];\dbS(\dbH))}\=\sup_{t\in
[t_1,t_2]}|F(t)|_{\cL(\dbH)}.
$$

Denote by $C_\cS([t_1,t_2];\dbS(\dbH))$ the set
of all strongly continuous mappings
$F:[t_1,t_2]\to \dbS(\dbH)$, that is,
$F(\cd)\xi$ is continuous on $[t_1,t_2]$ for
each $\xi\in H$. Let
$\{F_n\}_{n=1}^\infty\subset
C_\cS([t_1,t_2];\dbS(\dbH))$. We say
$\{F_n\}_{n=1}^\infty$ converges strongly to
$F\in C_\cS([t_1,t_2];\dbS(\dbH))$ if
$$
\lim_{n\to\infty}F_n(\cd)\xi=F(\cd)\xi,\qq
\forall \xi\in \dbH.
$$
In this case, we write
$$
\lim_{n\to\infty}F_n =F \q\mbox{ in
}\;C_\cS([t_1,t_2];\dbS(\dbH)).
$$
If $F\in C_\cS([t_1,t_2];\dbS(\dbH))$, then, by
the Uniform Boundedness theorem, the quantity
$$
|F|_{C_\cS([t_1,t_2];\dbS(\dbH))}\= \sup_{t\in
[t_1,t_2]}|F(t)|_{\cL(\dbH)}
$$
is finite, and $C_\cS([t_1,t_2];\dbS(\dbH))$ is
a Banach space with this norm (see
\cite{Bensoussan 1993} for the proof).

\ms

Let $H$ and $U$ be two separable Hilbert spaces.
Consider the following controlled linear
stochastic evolution equation (SEE, for short):
\begin{equation}\label{state}
\left\{
\begin{array}{ll}
\ns\ds dx =\big[(A+A_1 )x +B u
\big]ds+\big(C x +D u \big)dW(s)&\mbox{ in }(t,T], \\
\ns\ds x(t)=\eta\in
H,
\end{array}
\right.
\end{equation}
where  $A$ generates a $C_0$-semigroup
$\{e^{As}\}_{s\geq 0}$ on $H$,
$$\left\{\2n\begin{array}{ll}
\ns\ds A_1(\cd)\in L^1(0,T;\cL(H)),\q B(\cd)\in L^2(0,T;\cL(U;H)), \\
\ns\ds C(\cd)\in L^2(0,T;\cL(H)),\q D(\cd)\in
L^\i(0,T;\cL(H;U)).\end{array}\right.$$
In the above, $x(\cd)$ is the {\it state
process}, and $u(\cd)\in \cU[t,T]\=
L^2_\dbF(t,T;U)$  is the {\it control process}.
Any $u(\cd)\in\cU[t,T]$ is called an {\it
admissible control} (on $[t,T]$). For any {\it
initial pair} $(t,\eta)\in[0,T)\times H$ and
admissible control $u(\cd)\in\cU[t,T]$,
(\ref{state}) admits a unique mild solution
$x(\cd)\equiv x(\cd\,;t,\eta,u(\cd))$ (see Lemma
\ref{lm2}).  Here and in what follows, to
simplify the notations, the time variable $s$ is
suppressed in  $B$, $C$, etc.

Next we introduce the following cost functional:
\begin{equation}\label{cost}\begin{array}{ll}
\ns\ds \cJ(t,\eta;u(\cd))\deq\dbE \lan
Gx(T),x(T)\ran +\dbE \int_t^T \big(\lan Qx,x\ran
+ \lan Ru,u\ran \big) ds
,\end{array}\end{equation}
where
$$ G\in\dbS(H),\q Q(\cd)\in
L^1(0,T;\dbS(H)), \q R(\cd)\in
L^\infty(0,T;\dbS(U)). $$

The  optimal control problem studied in this
paper is as follows.

\ms

\bf Problem (SLQ). \rm For any given initial
pair $(t,\eta)\in[0,T)\times H$, find a $\bar
u(\cd)\in\cU[t,T]$, such that
\begin{equation}\label{optim}
V(t,\eta)\deq\cJ(t,\eta;\bar u(\cd))
=\inf_{u(\cd)\in\cU[t,T]}\cJ(t,\eta;u(\cd)).
\end{equation}

\medskip

Any $\bar u(\cd)\in\cU[t,T]$ satisfying
\eqref{optim} is called an {\it optimal control}
of Problem (SLQ) for the initial pair
$(t,\eta)$, and the corresponding $\bar
x(\cd)\equiv x(\cd\,; t,\eta,\bar u(\cd))$ is
called an {\it optimal state process}; the pair
$(\bar x(\cd),\bar u(\cd))$ is called an {\it
optimal pair}. The function $V(\cd\,,\cd)$ is
called the {\it value function} of Problem
(SLQ).

\begin{remark}
In this paper, we assume that $B(\cd)\in
L^2(0,T;\cL(U;H))$ and $D(\cd)\in L^\i(0,T;$
$\cL(H;U))$. Thus, our results can only be
applied to controlled stochastic partial
differential equations with distributed
controls. To study systems with boundary
controls, one needs to make some further
assumptions, such as the semigroup
$\{e^{As}\}_{s\geq 0}$ has some smoothing
effect. With such assumptions, we can prove
Theorem \ref{th4.6}(the details are lengthy and
beyond the scope of this paper). However, since
such assumptions contradict {\bf (AS2)} given
below, under them, we do not know how to prove
Theorem \ref{th4.4}. Fortunately, there are many
controlled stochastic partial differential
equations satisfying these conditions, such as
stochastic wave (\resp Schr\"odinger, KdV,
transport, plate) equations with internal
controls.
\end{remark}
\begin{remark}
In this paper, we assume that the coefficients
are deterministic. In such case, the
corresponding Riccati equation \eqref{Riccati}
is an operator-valued deterministic differential
equation. If one considers the problem that the
coefficients are stochastic, then an
operator-valued backward stochastic differential
equation should be studied. Until now, only some
very special cases of such equations are
investigated (e.g.\cite{GT 2005,GT 2014,LWZ}).
\end{remark}

The study of an optimal control problem for a
linear system with a quadratic cost functional
(LQ problem, for short) dates back at least to
\cite{Belman-Gicksberg-Gross 1958}, in which the
system is governed by a linear ordinary
differential equation. It an be regarded as the
simplest nontrivial optimal control problems,
namely, the system is linear and the cost
functional is quadratic. Consequently,  it has
elegant and fruitful mathematical structure.
Furthermore, it has  important applications
(e.g. \cite{Anderson-Moore 1989}). Such kind of
problem was investigated extensively in the
literature for a variety of deterministic
systems(e.g. \cite{Anderson-Moore 1989,Lasiecka
1991}).

LQ problems for controlled stochastic
differential equations (SDEs for short) was
first studied in \cite{Wonham 1968}. Such
problems are the most important examples of the
stochastic control problems, especially in their
applications in finance and economics. There are
a huge amount works addressing the LQ problems
for controlled SDEs(see \cite{Bismut
1976,Chen-Li-Zhou 1998, Chen-Zhou 2000, Hu-Zhou
2003, Kohlmann 2003, Qian-Zhou 2013, Rami 2001,
Tang 2003, Yong-Zhou 1999,Yu 2015} and the rich
references therein).

SEEs are used to describe a lot of random
phenomena appearing in physics, chemistry,
biology, and so on. In many situations SEEs are
more realistic mathematical models than the
deterministic ones (e.g.
\cite{Carmona,Kol,Kotelenez}). Thus, there are
many works addressing the optimal control
problems for SEEs. In particular, we refer the
readers to \cite{Ahmed 1981,Curtain 1977,Curtain
1978,GT 2005,GT 2014,Lasiecka 2017,Tessitore
1992,Tessitore 1992-1} and the rich references
therein for LQ problem for controlled SEEs.

In those works, the following assumption  was
taken for granted:
\begin{equation}\label{classical}
G\geq0,\qq R(s)\geq\d I,\qq Q(s)
\geq0,\qq\ae~s\in[t,T],
\end{equation}
for some $\d>0$. Under \eqref{classical}, when
all the operators in the SEEs are deterministic,
people had proven the corresponding Riccati
equation is uniquely solvable and Problem (SLQ)
admits a unique optimal control which has a
linear state feedback representation (under
certain technical
conditions)(e.g.\cite{Tessitore 1992}). On the
other hand, in the case that some operators in
the SEEs are stochastic, Problem (SLQ) was well
studied when $D=0$ and \eqref{classical} holds
(e.g. \cite{GT 2005,GT 2014}).

In \cite{Chen-Li-Zhou 1998}, the authors
discovered a new phenomenon, that is, Problem
(SLQ) might still be solvable for controlled
SDEs even if $R(s)$ is not positive definite for
a.e. $s\in[0,T]$. This motivated many subsequent
works concerning controlled SDEs (e.g.
\cite{Chen-Zhou 2000,Hu-Zhou 2003,Li
2002,Qian-Zhou 2013,Sun-Yong 2014,Sun-Li-Yong
2016,Yong-Zhou 1999}).  As far as we know, there
is no generalization of this to controlled SEEs,
which is one of the main purpose of this paper.

Recently, in \cite{Sun-Li-Yong 2016}, the
authors introduced the notions of open-loop and
closed-loop solvability of stochastic LQ
problems and showed the difference between these
two concepts. Roughly speaking, open-loop (\resp
closed-loop) solvability of a stochastic LQ
problem means that there is an open-loop (\resp
closed-loop) optimal control  of that problem.
Another main purpose of this paper is to study
the relationship between the closed-loop
solvability for Problem (SLQ) and the unique
solvability of the corresponding Riccati
equation.

In view of the main novel contributions
distinguishing this work from other publications
in the literatures are: \emph{(1) $R(s)\geq \d
I$ may not hold; (2) the equivalence between the
existence of an optimal feedback operator and
the existence of regular solution to the
stochastic Riccati equation is established.}
Although these two phenomena has been already
discovered for stochastic LQ problem for
controlled SDEs, one cannot simply mimic the
method to solve our problem. There are some
difficulties needing to be overcome. For
example, in finite dimensional case, one can
represent the solution to the Riccati equation
by the product of a solution to a matrix-valued
backward SDE and a solution to a matrix-valued
SDE. \emph{In the infinite-dimensional case,
formally, the matrix-valued SDE becomes an
operator-valued stochastic differential
equation.  Since there is no suitable
integration theory for general operator-valued
stochastic processes with respect to $W(\cd)$
(see \cite{vanNeerven1, vanNeerven2} for the
details), the solution to these operator-valued
processes cannot be defined in the classical
sense.}  More details on the difficulties
arising in the infinite dimensional settings can
be found in Sections
\ref{sec-pre}--\ref{sec-pr-main2}.

The rest of the paper is organized as follows.
In Section \ref{sec-main}, we present the main
results of this paper.   Section \ref{sec-pre}
is devoted to giving some preliminary results.
Sections \ref{sec-pr-main1} and
\ref{sec-pr-main2} are addressed to the proofs
of the main results. In Section
\ref{sec-convex}, we discuss the the uniform
convexity of the cost functional. At last, in
the appendix, we give proofs for some
preliminary results in Section \ref{sec-pre}.


\section{Statements of the main
results}\label{sec-main}


In this section, we present the main results of
this paper. To begin with, let us first
introduce some concepts.

\begin{definition}\label{def2.3}
We call $\cl\Th(\cd) \in L^2(t,T;\cL(U;H)) $ an
{\it optimal feedback operator} of Problem (SLQ)
on $[t,T]$ if
\begin{equation}\label{closed-opti}
\cJ(t,\eta;\cl\Th(\cd)\bar x(\cd))\leq
\cJ(t,\eta;u(\cd)),\qq\forall \eta\in H,\q
u(\cd)\in\cU[t,T],
\end{equation}
where $\bar x(\cd)$ is the mild solution to
\eqref{state} with $u(\cd)= \cl\Th(\cd)\bar
x(\cd)$.
\end{definition}
\begin{definition}\label{def2.3}
Problem (SLQ) is said to be ({\it uniquely})
{\it closed-loop solvable on $[t,T]$} if an
optimal feedback operator (uniquely) exists on
$[t,T]$.

Problem (SLQ) is said to be ({\it uniquely})
{\it closed-loop solvable} if it is (uniquely)
closed-loop solvable on any time horizon
$[t,T]$.
\end{definition}
\begin{remark}\label{rmk2.1}
Clearly, if $\cl\Th(\cd)$ is an optimal feedback
operator of Problem (SLQ) on $[t,T]$, then the
control $\bar u(\cd)\equiv\cl\Th(\cd)\bar
x(\cd)$ is an optimal control of Problem (SLQ)
for the initial pair $(t,\bar x(t))$. Hence, the
closed-loop solvability of Problem (SLQ) implies
the existence of an optimal control of that
problem. However, the converse is untrue. An
counterexample is given in \cite{Sun-Li-Yong
2016} for controlled SDEs.
\end{remark}

Next, let us recall the definition of the
generalized pseudo inverse of a self-adjoint
operator. More details and related proofs can be
found in \cite{Beutler 1965}.
\begin{definition}
Let $H_1$ be a Hilbert space. For $F\in
\cS(H_1)$, a generalized pseudo inverse of $F$
is defined as a linear operator $F^{\dag}:
\cD(F^{\dag}) \to H_1$ satisfying the following
four criteria:
$$
FF^{\dag}F=F,\q F^{\dag}FF^{\dag}=F^{\dag},\q
(FF^{\dag})^{*}=FF^{\dag},
\q(F^{\dag}F)^{*}=F^{\dag}F.
$$
\end{definition}
When $F$ is injective, $F^{\dag}$ is a left
inverse of $F$. In this case, $F^{\dag}F=I$.
When $F$ is surjective, $F^{\dag}$ is a right
inverse, as $FF^{\dag}=I$. When $F$ is
self-adjoint, $F^{\dag}$ always exists, which
may be unbounded.

Next, we introduce the Riccati equation
associated with Problem (SLQ) below:
\begin{equation}\label{Riccati}
\left\{\2n
\begin{array}{ll}
\ns\ds\dot P +P \big(A+A_1 \big)+\big(A+A_1
\big)^* P +C^* PC
+Q-L^* K^\dag L=0 &\mbox{ in }[t,T),\\
\ns\ds P(T)=G,
\end{array}
\right.
\end{equation}
where
$$
L(\cd)=B(\cd)^*P(\cd)+D(\cd)^* P(\cd)C(\cd),\qq
K(\cd)=R(\cd)+D(\cd)^* P(\cd)D(\cd).
$$
\begin{definition}\label{def3}
We call $P\in C_\cS([t,T];\dbS(H))$ a mild
solution to \eqref{Riccati} if for any $\eta\in
H$ and $s\in [t,T]$,
\begin{equation*}\label{8.20-eq23}
\begin{array}{ll}\ds
P(s)\eta\3n&\ds=e^{A^*(T-s)}Ge^{A(T-s)}\eta +
\int_s^T e^{A^*(\tau-s)}\big(P A_1\! +\! A_1^*
P\!+\! C^* PC \!+\!Q \!-\! L^*K^\dag
L\big)e^{A(\tau-s)}\eta d\tau.
\end{array}
\end{equation*}
\end{definition}

Now we can give the following notion.
\begin{definition}\label{def5}
A mild solution $P(\cd)\in C([t,T];\dbS(H))$ of
\eqref{Riccati} is {\it regular} if
\begin{equation}\label{regular-1}
\cR\big(L(s)\big)\subseteq\cR\big(K(s)\big),\qq
\ae~s\in[t,T],
\end{equation}
\begin{equation}\label{regular-2}
K(\cd)^\dag L(\cd) \in L^2(t,T;\cL(H;U)),
\end{equation}
and
\begin{equation}\label{regular-3}
K(s)\geq0,\qq\ae~s\in[t,T].
\end{equation}
\end{definition}
\begin{definition}\label{def6}
A solution $P(\cd)$ of \eqref{Riccati} is {\it
strongly regular} if
\begin{equation}\label{strong-regular}
K(s)\geq \l I,\qq\ae~s\in[t,T],
\end{equation}
for some $\l>0$.
\end{definition}
\begin{definition}\label{def6.1}
The Riccati equation (\ref{Riccati}) is said to
be ({\it strongly }) {\it regularly solvable},
if it admits a (strongly) regular solution.
\end{definition}
\begin{remark}
As far as we know, the notion of (strongly)
regular solution of a matrix-valued Riccati
equation was first introduced in
\cite{Sun-Li-Yong 2016}.
\end{remark}
\begin{remark}
Clearly, condition (\ref{strong-regular})
implies (\ref{regular-1})--(\ref{regular-3}).
Thus, a strongly regular solution $P(\cd)$ is
regular. Moreover, if $P$ is a strongly regular
solution, then $K(s)$ is invertible for a.e.
$s\in[0,T]$, i.e., the generalized pseudo
inverse is the inverse of $K(s)$.
\end{remark}

To investigate the relation between the
closed-loop solvability of Problem {\rm(SLQ)}
and the existence of a regular solution to the
Riccati equation {\rm(\ref{Riccati})}, we need
to make the following assumptions for $A$.

\vspace{0.1cm}

{\bf (AS1)} $A$ generates a $C_0$-group on $H$.

\vspace{0.1cm}

{\bf (AS2)} The eigenfunctions
$\{e_j\}_{j=1}^\infty$ of $A$ constitutes an
orthonormal basis of $H$.

\vspace{0.1cm}

\begin{remark}
We put two assumptions on $A$, i.e., $A$
generates a $C_0$-group and its eigenfunctions
constitutes an orthonormal basis of $H$. Both of
them play important roles in the proof of
Theorem \ref{th4.4}. Indeed, in the proof of
Theorem \ref{th4.4}, we need to find a finite
dimensional approximation of \eqref{Riccati},
i.e., we should approximate operator-valued
processes by matrix-valued processes.  To this
end,   $H$ should has an orthonormal basis
$\{e_j\}_{j=1}^\infty$. This is true since $H$
is separable. However, for getting some good
estimates (see Steps 1 and 2 in the proof of
Theorem \ref{th4.4}), we need the fact that
$e_j$ is an eigenfunction of $A$ for each
$j\in\dbN$. This leads to {\bf (AS2)}
Furthermore, to get the inverse of the operator
$X(s,\cd)$ given in step 2 in the proof of
Theorem \ref{th4.4}, we need {\bf (AS1)}.

It seems that both {\bf (AS1)} and {\bf (AS2)}
are only technical assumptions and can be
dropped. However, we do not know how to do it
now. For example, without {\bf (AS1)}, $-A^*$
does not generate a $C_0$-semigroup. Then the
equation \eqref{5.26-eq4} is not well-posed and
the operator $\wt X(s,\cd)$ is not well defined.
Then one cannot show the inverse of $X(s,\cd)$
On the other hand, without {\bf (AS2)}, we
cannot get \eqref{6.8-eq5.1} and
\eqref{6.8-eq5.2}, which are  keys in the proof
of many results, such as \eqref{7.24-eq4} and
\eqref{7.24-eq3}. More details can be see in the
proof of Theorem \ref{th4.4}.

Fortunately, under these conditions, the system
\eqref{state} covers many controlled stochastic
PDEs, such as stochastic wave equations,
stochastic Schr\"odinger equations, stochastic
beam equations, with internal controls.
\end{remark}
\begin{theorem}\label{th4.4}

i) If the Riccati equation {\rm(\ref{Riccati})}
admits a regular solution $P(\cd)\in
C_\cS([t,T];\dbS(H))$, then Problem {\rm(SLQ)}
is closed-loop solvable. In this case, the
optimal feedback operator $\cl\Th(\cd)$ is given
by
\begin{equation}\label{Th-v-rep}
\cl\Th=-K^\dag L +\big(I-K^\dag L\big)\th,
\end{equation}
for some $\th(\cd)\in L^2(t,T;\cL(U;H))$, and
the value function is
\begin{equation}\label{Value}
V(t,\eta)=\dbE \langle P(t)\eta,\eta\rangle.
\end{equation}

ii) Let {\bf (AS1)} and {\bf (AS2)} hold. If
Problem {\rm(SLQ)} is closed-loop solvable, then
the Riccati equation {\rm(\ref{Riccati})} admits
a regular solution $P(\cd)\in
C_\cS([t,T];\dbS(H))$.
\end{theorem}

Next result gives a sufficient and necessary
condition for the existence of a strongly
regular solution to the Riccati equation
{\rm(\ref{Riccati})}.

\begin{theorem}\label{th4.6}
The following statements are equivalent:

\ms

{\rm(i)} The map $u(\cd)\mapsto \cJ(0,0;u(\cd))$
is uniformly convex, i.e., there exists a $\l>0$
such that
\begin{equation}\label{J>l*}
\cJ(0,0;u(\cd))\geq\l\,\dbE\1n\int_0^T|u(s)|_U^2ds,\qq\forall
u(\cd)\in\cU[0,T],
\end{equation}

\ms

{\rm(ii)} The Riccati equation
{\rm(\ref{Riccati})} admits a strongly regular
solution $P(\cd)\in C_\cS([0,T];\dbS(H))$.
\end{theorem}
\begin{remark}
Clearly, if \eqref{classical} holds, then the
map $u(\cd)\mapsto \cJ(0,0;u(\cd))$ is uniformly
convex. On the other hand, there are some
interesting cases that the map $u(\cd)\mapsto
\cJ(0,0;u(\cd))$ is uniformly convex but
\eqref{classical} does not hold. Please see
Section \ref{sec-convex} for the details.
\end{remark}

Combining Theorems \ref{th4.4} and  \ref{th4.6},
we obtain the following corollary.

\begin{corollary}\label{cor4.8}
Let \eqref{J>l*} hold. Then, at any
$(t,\eta)\in[0,T)\times H$, Problem {\rm(SLQ)}
admits a unique  optimal control $\bar u(\cd)$
of a state feedback form:
\begin{equation}\label{opti-biaoshi}
\bar u(\cd)=-K(\cd)^{-1}L(\cd)\bar x(\cd),
\end{equation}
where $P(\cd)$ is the unique strongly regular
solution of {\rm(\ref{Riccati})} with $\bar
x(\cd)$ being the solution to the following
closed-loop system:
\begin{equation}\label{closed-loop-state}
\left\{\2n
\begin{array}{ll}
\ds d\bar x = \big(A+A_1-BK^{-1}L\big)\bar xds +
\big(C-DK^{-1}L\big)\bar
x dW(s) &\mbox{ in }(t,T], \\
\ns\ds \bar x(t)=\eta.
\end{array}
\right.
\end{equation}
\end{corollary}

{\it Proof}\,: By Theorem \ref{th4.6}, the
Riccati equation (\ref{Riccati}) admits a unique
strongly regular solution $P(\cd)\in
C([0,T];\dbS(H))$.  Applying Theorem
\ref{th4.4}, we get the desired result.
\endpf


\section{Some preliminaries}\label{sec-pre}


In this section, we present some useful results
which will be used  in the sequel. Except Lemma
\ref{lm2.3}, the proofs of  other results are
put in the appendix.

\ms

For any $t\in[0,T)$, consider the following SEE:
\begin{equation}\label{6.20-eq1}
\left\{
\begin{array}{ll}\ds
dx = [(A+\cA) x  + f]ds + (\cB x+g)dW(s)
&\mbox{ in }(t,T],\\
\ns\ds x(t)=\eta.
\end{array}
\right.
\end{equation}
Here $\cA \in L^1(t,T;\cL(H))$, $\cB\in
L^2(t,T;\cL(H))$, $\eta\in L^2_{\cF_t}(\Om;H)$,
$f\in  L^2_\dbF(\Om;L^1(t,T;H))$  and $g\in
L^2_\dbF(t,T;H)$.

\begin{lemma}\label{lm2}
The equation \eqref{6.20-eq1} admits a unique
mild solution $x(\cd)\in
C_\dbF([t,T];L^2(\Om;H))$. Moreover,
\begin{equation}\label{lm2-eq1}
|x(\cd)|_{C_\dbF([t,T];L^2(\Om;H))}\leq
\cC\big(|\eta|_{L^2_{\cF_t}(\Om;H)} +
|f|_{L^2_\dbF(\Om;L^1(t,T;H))} +
|g|_{L^2_\dbF(t,T;H)}\big).
\end{equation}
\end{lemma}

Next, consider the following backward stochastic
evolution equation (BSEE for short):
\begin{equation}\label{6.20-eq10}
\left\{
\begin{array}{ll}\ds
dy = -[(A+A_1)^*y+ C^* z + h]ds + zdW(s) &\mbox{ in }[t,T),\\
\ns\ds y(T)=\xi.
\end{array}
\right.
\end{equation}
Here $\xi\in L^2_{\cF_T}(\Om;H)$  and $h\in
L^2_\dbF(\Om;L^1(0,T;H))$. We have the following
result:
\begin{lemma}\label{lm3}
The equation \eqref{6.20-eq10} admits a unique
mild solution $(y(\cd),z(\cd))\in
L^2_\dbF(\Om;C([0,T];$ $H))\times
L^2_\dbF(0,T;H)$, and
\begin{equation}\label{2.9-eq0}
|(y(\cd),z(\cd))|_{L^2_\dbF(\Om;C([0,T];H))\times
L^2_\dbF(0,T;H)}\leq
\cC\big(|\xi|_{L^2_{\cF_T}(\Om;H)}+|h|_{L^2_\dbF(\Om;L^1(0,T;H))}\big).
\end{equation}
\end{lemma}

Next, we recall the following result.

\begin{lemma}\label{5.7-prop1}
Let $\cl\Th(\cd)$ be an optimal feedback
operator of Problem (SLQ). Let $\zeta\in H$. The
following forward-backward stochastic evolution
equation (FBSEE for short) admits a mild
solution $(\bar x(\cd), \bar y(\cd),\bar
z(\cd))\in C_\dbF([t,T];L^2(\Om;H))\times
C_\dbF([t,T];L^2(\Om;H))\times L^2_\dbF(t,T;H)$:
\begin{equation}\label{FBSDE5.1}
\left\{
\begin{array}{ll}
\ns\ds d x =(A+A_1+B\cl\Th) x ds+ (C+D\cl\Th) x
dW(s)
&\hb{\rm in } (t,T],\\
\ns\ds
d y =-\big[(A+A_1)^* y+C^*  z+ Q x\big]ds+ zdW(s)&\mbox{\rm in } [t,T),\\
\ns\ds  x(t)=\zeta,\qq  y(T)=G x(T),
\end{array}
\right.
\end{equation}
and the following condition holds:
\begin{equation}\label{5.7-eq3}
B^* \bar y +D^* \bar z + R\Th \bar x =0,\q \ae
(s,\om)\in [t,T]\times\Om.
\end{equation}
\end{lemma}
If $A_1,C\in L^\infty(0,T;\cL(H))$ and $B,D\in
L^\infty(0,T;\cL(U;H))$, Lemma \ref{5.7-prop1}
is a trivial corollary of Theorem 5.2 in
\cite{LZ}. The general case can be handled
similarly. For the readers' convenience, we give
the proof in the appendix.

\vspace{0.1cm}

\vspace{0.1cm}

Consider the following SEE:
\begin{equation}\label{5.26-eq4}
\left\{
\begin{array}{ll}
\ns\ds d\tilde x =
\big[-A-A_1-B\cl\Th+\big(C+D\cl\Th\big)^2
\big]^*\tilde x ds - \big(C+D\cl\Th\big)^*\tilde
xdW(s) &  \mbox{ in }(t,T],\\
\ns\ds \tilde x(t)=\zeta,
\end{array}
\right.
\end{equation}
where $\zeta\in H$. If $A$ generates a
$C_0$-group, then $-A^*$ also generates a
$C_0$-group. In this case, by Lemma \ref{lm2},
the equation \eqref{5.26-eq4} admits a unique
mild solution $\tilde x(\cd)\in
C_\dbF([t,T];L^2(\Om;H))$.

Let $\{\f_i\}_{i=1}^\infty$ be an orthonormal
basis of $U$. For each $n\in\dbN$, denote by
$\G_n$ (\resp $\wt \G_n$) the projection
operator from $H$ (\resp $U$) to
$H_n\=\span_{1\leq j\leq n}\{e_j\}$ (\resp
$U_n\=\span_{1\leq j\leq n}\{\f_j\}$). Write
\begin{equation}\label{8.20-eq8}
\begin{array}{ll}\ds
A_{n}=\G_nA\G_n,\q A_{1,n}=\G_nA_1\G_n,\q B_n =
\G_n B \wt \G_n,\q C_n
= \G_n C \G_n, \;\;\; D_n = \G_n D \wt\G_n,\\
\ns\ds  Q_n = \G_n Q \G_n,\;\;\; R_n = \wt\G_n R
\wt\G_n,\;\;\; G_n=\G_nG\G_n,\;\;\;\Th_n =
\wt\G_n \cl\Th \G_n.
\end{array}
\end{equation}
Denote by $\{\l_j\}_{j=1}^\infty$ the
eigenvalues of $A$ such that $Ae_j=\l_j e_j$
(recall that $\{e_j\}_{j=1}^\infty$ is the
eigenfunctions of $A$). Then,
\begin{equation}\label{6.8-eq5.2.1}
e^{A_ns}\zeta = \sum_{j=1}^n e^{\l_js}\lan
\zeta, e_j\ran e_j = e^{As}\G_n\zeta,\qq \forall
\zeta\in H.
\end{equation}
Thus,
\begin{equation}\label{6.8-eq5.1}
\lim_{n\to+\infty} e^{A_{n}s}\zeta = e^{As}
\zeta,\qq \forall\, s\in [0,T],\q\zeta\in H.
\end{equation}
Similarly, we can get that
\begin{equation}\label{6.8-eq5.2}
\lim_{n\to+\infty} e^{A^*_{n}s}\zeta = e^{A^*s}
\zeta,\qq \forall\, s\in [0,T],\q\zeta\in H.
\end{equation}
Further, it is easy to show that
\begin{equation}\label{6.8-eq5}
\begin{array}{ll}\ds
\lim_{n\to+\infty} A_{1,n}\zeta = A_1 \zeta, \q
\lim_{n\to+\infty} C_n\zeta = C \zeta, \q
\lim_{n\to+\infty} G_n\zeta = G \zeta, \\
\ns\ds \lim_{n\to+\infty} Q_n\zeta = Q \zeta,
\qq
\lim_{n\to+\infty} \Th_n\zeta = \cl\Th \zeta,\\
\ns\ds \hspace{4.7cm}\hb{ for all }\zeta\in H
\hb{ and }\ae (s,\om)\in [0,T]\times\Om,
\end{array}
\end{equation}
\begin{equation}\label{6.8-eq5.1.1}
\begin{array}{ll}\ds
\lim_{n\to+\infty} B_n\varsigma = B \varsigma,
\q \lim_{n\to+\infty} D_n\varsigma = D
\varsigma,\q\lim_{n\to+\infty} R_n\varsigma= R \varsigma, \\
\ns\ds \hspace{6cm}\hb{ for all }\varsigma\in U
\hb{ and }\ae (s,\om)\in [0,T]\times\Om.
\end{array}
\end{equation}

Let $\zeta\in H$, consider the following
equations:
\begin{equation}\label{6.7-eq1}
\left\{
\begin{array}{ll}
\ds dx_n=  (A_n+A_{1,n}+B_n\Th_n) x_n ds+
(C_n+D_n\Th_n) x_n
dW(s)& \mbox{ in }[t,T],\\
\ns\ds dy_n=-\big[(A_n+A_{1,n})^* y_n+ C^*_n
z_n+ Q_nx_n\big]ds+ z_ndW(s) &
\mbox{ in }[t,T],\\
\ns\ds x_n(t)=\G_n\zeta, \q y_n(T)= G_n x_n(T)
\end{array}
\right.
\end{equation}
and
\begin{equation}\label{6.7-eq2}
\left\{\!\!\!
\begin{array}{ll}
\ds d\tilde x_n=\!
\big[\!-\!A_n\!\!-\!A_{1,n}\!-\!B_n\Th_n+\big(C_n\!+\!D_n\Th_n\big)^2
\big]^*\tilde x_n ds \!-
(C_n\!+\!D_n\Th_n)^*\tilde x_n
dW(s)& \mbox{in }(t,T],\\
\ns\ds \tilde x_n(0)=\G_n\zeta.
\end{array}
\right.
\end{equation}

We have the following result.
\begin{lemma}\label{lm4}
Let {\bf (AS1)} and {\bf (AS2)} hold. For any
$\zeta\in H$, it holds that
\begin{equation}\label{6.7-eq3}
\left\{
\begin{array}{ll}\ds
\lim_{n\to+\infty}  x_n(\cd) = x(\cd) &\mbox{ in
}
C_\dbF([t,T];L^2(\Om;H)),\\
\ns\ds \lim_{n\to+\infty} y_n(\cd) = y(\cd)
&\mbox{ in }
L^2_\dbF(\Om;C([t,T];H)), \\
\ns\ds \lim_{n\to+\infty} z_n(\cd) = z(\cd)
&\mbox{ in } L^2_\dbF(t,T;H),
\\
\ns\ds \lim_{n\to+\infty} \tilde x_n(\cd)
=\tilde x(\cd) &\mbox{ in }
C_\dbF([t,T];L^2(\Om;H)).
\end{array}
\right.
\end{equation}
\end{lemma}

\vspace{0.41cm}

Let $\Th(\cd)\in L^2(t,T;\cL(H;U))$. Consider
the following operator-valued equation:
\begin{equation}\label{eq-Lya}
\left\{\2n\begin{array}{ll}
\ds\dot P+P(A+A_1+B\Th)+(A+A_1+B\Th)^* P\\
\ns\ds+(C+D\Th)^* P(C+D\Th)+\,\Th^* R\Th +Q=0 &\mbox{ in }[t,T),\\
\ns\ds P(T)=G.
\end{array}
\right.
\end{equation}
\begin{definition}\label{def1}
We call $P\in C_\cS([t,T];\dbS(H))$ a mild
solution to \eqref{eq-Lya} if for any $s\in
[t,T]$,
\begin{equation}\label{8.20-eq22}
\begin{array}{ll}\ds
P(s)\eta\3n&\ds=e^{(T-s)A^*}Ge^{(T-s)A}\eta +
\int_s^T
e^{(\tau-s)A^*}[P(A_1+B\Th)+(A_1+B\Th)^*
P\\
\ns&\ds\q+(C+D\Th)^* P(C+D\Th)+\,\Th^* R\Th
+Q]e^{(\tau-s)A}\eta d\tau,\qq \forall\, \eta\in
H.
\end{array}
\end{equation}
\end{definition}
\begin{proposition}\label{prop3}
There is a unique mild solution to
\eqref{eq-Lya}. Moreover,
\begin{equation*}\label{prop3-eq1}
\!\!\begin{array}{ll}\ds
|P|_{C_\cS([t,T];\cS(H))} \!\leq \!\cC
e^{\int_t^T (2|A_1\!+\!B\Th|_{\cL(H)}
+|C\!+\!D\Th|_{\cL(H)}^2 ) ds}\[|G|_{\cL(H)}
\!+\! \int_t^T \!\!\big( |\Th|_{\cL(H;U)}^2
|R|_{\cL(U)} \!+\!|Q|_{\cL(H)}\big)ds\].
\end{array}
\end{equation*}
\end{proposition}

The following result illustrates the
differentiability of $P$.
\begin{proposition}\label{prop1}
Let $P$ be a mild solution to \eqref{eq-Lya}.
Then for any $\eta,\zeta\in D(A)$, $\lan
P(\cd)\eta, \zeta\ran$ is differentiable in
$[t,T]$ and
\begin{equation}\label{8.20-eq20}
\begin{array}{ll}\ds
\frac{d}{ds}\lan P\eta, \zeta\ran\3n&\ds=-\lan
P\eta, (A+A_1+B\Th)\zeta\ran - \lan
P(A+A_1+B\Th)\eta, \zeta\ran \\
\ns&\ds\q - \lan P(C+D\Th)\eta, (C+D\Th)\zeta
\ran- \lan R\Th \eta,  \Th \zeta \ran- \lan Q
\eta, \zeta \ran.
\end{array}
\end{equation}
\end{proposition}

Similarly, we have the following result.
\begin{proposition}\label{prop2}
Let $P$ be a mild solution to \eqref{Riccati}.
Then for any $\eta,\zeta\in D(A)$, $\lan
P(\cd)\eta, \zeta\ran$ is differentiable in
$[t,T]$ and
\begin{equation}\label{8.20-eq24}
\begin{array}{ll}\ds
\frac{d}{ds}\lan P\eta, \zeta\ran\3n&\ds=-\lan
P\eta, (A+A_1)\zeta\ran - \lan
P(A+A_1)\eta, \zeta\ran \\
\ns&\ds\q - \lan P C \eta, C \zeta \ran- \lan Q
\eta, \zeta \ran+\lan K^\dag L \eta, L \zeta
\ran.
\end{array}
\end{equation}
\end{proposition}

The  result below gives a relation between the
cost functional and the Riccati equation
\eqref{Riccati}.

\begin{lemma}\label{lm2.3}
Let $P(\cd)\in C_\cS([0,T];\dbS(H))$ be the mild
solution to \eqref{eq-Lya} with $t=0$. Then for
any $(t,\eta)\in[0,T)\times H$ and
$u(\cd)\in\cU[t,T]$, we have
\begin{equation}\label{8.20-eq29}
\begin{array}{ll}
\ns\ds \cJ(t,\eta;\Th(\cd)x(\cd)+u(\cd))
=\langle P(t)\eta,\eta\rangle
+\dbE\int_t^T\big[2\lan\big(L+K\Th\big)x,u\ran
+\lan Ku,u\ran \big]ds.
\end{array}
\end{equation}
\end{lemma}

In the proof of Lemma \ref{lm2.3}, we use a
density argument, i.e., approximating the
solution $x$ by a $D(A)$-valued process $x_\l$,
doing calculation for $x_\l$ and letting $\l$
tend to $+\infty$ to get the equality for $x$.
The main reason is that one may not apply
It\^o's formula to $|x(t)|_H^2$ and $x(t)$ may
not take values in $D(A)$. This technique will
be used several times in the rest of this
papers. Except the proof for Lemma \ref{lm2.3},
we omit such process and apply It\^o's formula
to $|x(t)|_H^2$ directly and assume that $x(t)$
belongs to $D(A)$.  Hence, we give the proof
here rather than put it in the appendix.

\vspace{0.2cm}

{\it Proof of Lemma \ref{lm2.3}}: For any
$(t,\eta)\in[0,T)\times H$ and
$u(\cd)\in\cU[t,T]$, let $x(\cd)$ be the
solution to
$$
\left\{
\begin{array}{ll}
\ns\ds dx =\big[(A+A_1+B\Th)x+Bu\big]ds
+\big[(C+D\Th)x+Du\big]dW(s) &\mbox{ in }(t,T], \\
\ns\ds x(t)=\eta.
\end{array}
\right.
$$
Let $\dbR(\l)\= \l I(\l I-A)^{-1}$ for
$\l\in\rho(A)$ (the resolvent of $A$) and
$x_\l(\cd)=\dbR(\l) x(\cd)$. Then $x_\l(\cd)$ is
the mild solution to
$$
\left\{\2n
\begin{array}{ll}
\ds dx_\l =\big\{Ax_\l \!+\!\dbR(\l)\big[(A_1\!+\!B\Th)x\!+\!Bu\big]\big\}ds \!+\!\dbR(\l)\big[(C\!+\!D\Th)x\!+\!Du\big]dW(s) &\mbox{ in }(t,T], \\
\ns\ds x_\l(t)=\dbR(\l)\eta.
\end{array}
\right.
$$

By It\^o's formula and Proposition \ref{prop1},
we have
\begin{eqnarray}\label{8.20-eq28}
&& \q\mE\[\lan Gx_\l(T),x_\l(T)\ran
+\int_t^T\big(\lan Qx_\l,x_\l\ran + \lan R(\Th
x_\l+u),\Th x_\l+u\ran\big)
ds\]\nonumber\\
&&\ds=\!\lan P(t)\dbR(\l)\eta,\dbR(\l)\eta\ran
\!+\!\dbE\!\int_t^T\!\!\!\Big\{\!\!-\!\lan
\!P(t)x_\l, (A\!+\!A_1\!\!+\!\!B\Th)x_\l\ran \!-\! \lan P(s)(A\!+\!A_1\!\!+\!\!B\Th)x_\l, x_\l\ran\nonumber\\
&&\ds\q   - \lan P(C+D\Th)x_\l, (C+D\Th)x_\l
\ran-
\lan R\Th x_\l,  \Th x_\l \ran- \lan Q x_\l, x_\l\ran +\lan P A x_\l,x_\l\ran\nonumber\\
&& \ds\q+\lan
P\dbR(\l)\big[(A_1+B\Th)x+Bu\big],x_\l\ran+\lan
Px_\l, Ax_\l\ran
+\lan Px_\l,\dbR(\l)[(A+B\Th)x+Bu]\ran\nonumber\\
&&\ds\q  +\lan P\dbR(\l)\big[(C+D\Th)x+Du\big],\dbR(\l)\big[(C+D\Th)x+Du\big]\ran+\lan Qx_\l,x_\l\ran\nonumber\\
&&\ds\q +\lan R(\Th x_\l+u),\Th x_\l+u\ran\Big\}ds \\
&&\ds=\lan P(t)\eta,\eta\ran
+\dbE\int_t^T\big[2\lan\big(L+K\Th\big)x,u\ran
+\lan Ku,u\ran\big]ds + F(\l),\nonumber
\end{eqnarray}
where
\begin{eqnarray*}
&&\ds \q F(\l) =\langle
P(t)\dbR(\l)\eta,\dbR(\l)\eta\rangle-\langle
P(t)\eta,\eta\rangle +\!\dbE\!\int_t^T\!\!\Big\{
\lan P\dbR(\l)(A_1+B\Th)x,x_\l\ran\\
&&\ds\q -\lan P (A_1+B\Th)x_\l,x_\l\ran
+\lan Px_\l,\dbR(\l) (A+B\Th)x \ran - \lan P(s)(A_1\!+\!B\Th)x_\l, x_\l\ran\\
&&\ds\q  +\lan
P\dbR(\l)(C+D\Th)x,\dbR(\l)(C+D\Th)x\ran- \lan
P(C+D\Th)x_\l, (C+D\Th)x_\l \ran \\
&&\ds\q + \lan P\dbR(\l) Bu,x_\l\ran \!-\! \lan
PBu,x\ran \! +\!\lan Px_\l,\dbR(\l) Bu \ran
\!-\! \lan
Px,Bu \ran \!+\!\lan P\dbR(\l)Du,\dbR(\l)Du\ran\\
&&\ds\q - \lan P Du, Du\ran +\lan
P\dbR(\l)(C+D\Th)x,\dbR(\l)Du\ran - \lan
P(C+D\Th)x,Du\ran\\
&&\ds\q +\lan P\dbR(\l)Du,\dbR(\l)(C+D\Th)x\ran
- \lan PDu,(C+D\Th)x\ran\Big\}ds.
\end{eqnarray*}
Noting that for any $\zeta\in H$,
\begin{equation}\label{8.20-eq26.1}
\lim_{\l\to\infty}\dbR(\l)\zeta = \zeta \q\mbox{
in }H,
\end{equation}
we have that
$$
\lim_{\l\to\infty}\langle
P(t)\dbR(\l)\eta,\dbR(\l)\eta\rangle = \langle
P(t)\eta,\eta\rangle
$$
and
\begin{equation}\label{8.20-eq26}
\lim_{\l\to\infty} x_\l= x \q\mbox{ in
}\;C_\dbF([t,T];L^2(\Om;H)).
\end{equation}
By \eqref{8.20-eq26.1} and \eqref{8.20-eq26}, we
get that for a.e. $s\in [t,T]$,
\begin{equation}\label{8.20-eq27}
\begin{array}{ll}\ds
\lim_{\l\to\infty}\big[\lan
P(s)\dbR(\l)(A_1(s)+B(s)\Th(s))x(s),x_\l(s)\ran\\
\ns\ds\qq -\lan P
(s)(A_1(s)+B(s)\Th(s))x_\l(s),x_\l(s)\ran\big]=0,\q\dbP\mbox{-a.s.}
\end{array}
\end{equation}
It follows from the definition of $\dbR(\l)$
that
$$
\begin{array}{ll}\ds
\big|\lan
P(s)\dbR(\l)(A_1(s)+B(s)\Th(s))x(s),x_\l(s)\ran
-\lan P
(s)(A_1(s)+B(s)\Th(s))x_\l(s),x_\l(s)\ran\big|\\
\ns\ds\leq \cC |P(s)|_{\cL(H)}
(|A_1(s)|_{\cL(H)}+|B(s)|_{\cL(U;H)}|\Th(s)|_{\cL(H;U)})|x(s)|_H^2.
\end{array}
$$
This, together with \eqref{8.20-eq27} and
Lebesgue's dominated convergence theorem,
implies that
$$
\lim_{\l\to\infty}\dbE \int_t^T\big[\lan
P\dbR(\l)(A_1+B\Th)x,x_\l\ran -\lan P
(A_1+B\Th)x_\l,x_\l\ran \big] ds = 0.
$$
By a similar argument, we can prove that $
\ds\lim_{\l\to+\infty}F(\l) = 0$.
Letting $\l$ tend to $+\infty$ in both sides of
\eqref{8.20-eq28}, we get \eqref{8.20-eq29}.
\endpf

\begin{remark}
In the derivation of \eqref{8.20-eq28}, we use
the fact that the mean value of a stochastic
integral of a function quadratically depending
on $x$ is zero.
\end{remark}
\begin{remark}
Since $x(\cd)$ may not be $D(A)$-valued, in the
proof of Lemma \ref{lm2.3}, we introduce a
family of $\{x_\l\}_{\l\in\rho(A)}$ to apply
Proposition \ref{prop1}. In the rest of this
paper, we omit such procedures to save the space
and simply apply Proposition \ref{prop1} to
$\lan P(\cd)x(\cd),x(\cd)\ran$.
\end{remark}
\vspace{0.3cm}

Next, we give a result concerning the existence
of an optimal control of Problem {\rm(SLQ)}.
\begin{proposition}\label{prop4.1}
Suppose the map $u(\cd)\mapsto \cJ(0,0;u(\cd))$
is uniformly convex. Then Problem {\rm(SLQ)}
admits a unique optimal control, and there
exists a constant $\a\in\dbR$ such that
\begin{equation}\label{uni-convex-prop0}
V(t,\eta)\geq\a|\eta|^2,\qq\forall
(t,\eta)\in[0,T]\times H.
\end{equation}
\end{proposition}

The next result shows that the solution to
\eqref{eq-Lya} is bounded below.
\begin{proposition}\label{prop4.5}
Let  {\rm(\ref{J>l*})} hold. Then for any
$\Th(\cd)\in L^2(t,T;\cL(U;H))$, the solution
$P(\cd)\in C_\cS([t,T];\dbS(H))$ to
\eqref{eq-Lya} satisfies
\begin{equation}\label{Convex-prop-1}
K(s)\geq\l I, \q\ae~s\in[t,T],\q \hb{and}\q
P(s)\geq \a I,\q\forall s\in[t,T],\end{equation}
where $\a\in\dbR$ is the constant appearing in
{\rm(\ref{uni-convex-prop0})}.
\end{proposition}

Let
$$\wt A(\cd)\in L^1(0,T;\cL(H)),\q \wt C(\cd)\in L^2(0,T;\cL(H)),
\q\wt Q(\cd)\in L^1(0,T;\dbS(H)), \q \wt
G\in\dbS(H).$$
Consider the following Lyapunov equation:
\begin{equation}\label{lm2.4-eq1}
\left\{ \2n\begin{array}{ll}
\ns\ds \dot{P} +P (A+\wt A )+(A+\wt A )^* P +\wt C^* P\wt C+\wt Q=0 &\mbox{ in }[t,T),\\
\ns\ds P(T)=\wt G.
\end{array}
\right.
\end{equation}

\begin{lemma}\label{lm2.4}
The equation \eqref{lm2.4-eq1} admits a unique
solution $P(\cd) \in C_\cS([t,T];\dbS(H))$.
Moreover, if
\begin{equation}\label{lm2.4-eq2}
\wt G\geq0,\qq\wt Q(s)\geq0,\qq \ae~s\in[t,T],
\end{equation}
then $P(\cd)\in C_\cS([t,T];\cl{\dbS_+}(H))$.
\end{lemma}
\begin{lemma}\label{lm2.5}
For any $u(\cd)\in\cU[t,T]$, let $x$ be the
corresponding solution to \eqref{state} with
$\eta = 0$.  Then for every $\Th(\cd)\in
L^2(t,T;\cL(H;U))$, there exists a constant
$c_0>0$ such that
\begin{equation}\label{lem-2.6}
\dbE \int_t^T\big|u(s)-\Th(s) x(s)\big|^2ds \geq
c_0\dbE\int_t^T|u(s)|^2ds,\qq\forall
u(\cd)\in\cU[t,T].
\end{equation}
\end{lemma}
%


\section{Proof of Theorem
\ref{th4.4}}\label{sec-pr-main1}


This section is devoted to the proof of Theorem
\ref{th4.4}.

\ms

\it Proof of Theorem \ref{th4.4}. \rm Proof of
assertion i). Take any $u(\cd)\in\cU[t,T]$, let
$x(\cd)\equiv x(\cd\,;t,\eta,u(\cd))$ be the
corresponding state process. Then
$$\begin{array}{ll}
\ns\ds
\cJ(t,\eta;u(\cd))=\frac{1}{2}\dbE\Big[\langle
Gx(T),x(T)\rangle  +\int_t^T\big(\big\langle
Qx,x\big\rangle
+\big\langle Ru,u\big\rangle \big)ds\Big]\\
\ns\ds=\frac{1}{2}\dbE \big\langle
P(t)\eta,\eta\big\rangle\!+\frac{1}{2}\dbE
\int_t^T \!\big[\big\langle\big(\!-P(A+A_1)\!
-\! (A+A_1)^* P \!-\! C^*
PC \!-\! Q \!+\! L^*K^\dag L\big)x,x\big\rangle\\
\ns\ds\qq +\big\langle P((A+A_1)x +
Bu),x\big\rangle+ \big\langle Px,(A+A_1)x +
Bu\big\rangle + \big\langle P(Cx + Du),Cx+ Du
\big\rangle\\
\ns\ds\qq +\big\langle
Qx,x\big\rangle+\big\langle
Ru,u\big\rangle\big]ds\\
\ns\ds=\frac{1}{2}\dbE\Big[\big\langle
P(t)\eta,\eta\big\rangle + \int_t^T\!\big(
\big\langle L^*K^\dag Lx,x\big\rangle +
2\big\langle Lx,u\big\rangle+\big\langle K
u,u\big\rangle\big)ds\Big].
\end{array}
$$
Noting that $L= B^* P+D^* PC=-K\Th$,
it holds that
$$
\begin{array}{ll}\ds
\q \cJ(t,\eta;u(\cd))\\
\ns\ds=\frac{1}{2}\dbE\Big[\big\langle
P(t)\eta,\eta\big\rangle
+\int_t^T\big(\big\langle L^* K^\dag
Lx,x\big\rangle+2\big\langle Lx,u\big\rangle
+\big\langle Ku,u\big\rangle\big)ds\Big]\\
\ns\ds=\frac{1}{2}\dbE\Big[\big\langle
P(t)\eta,\eta\big\rangle+\int_t^T\big(
\big\langle \Th^* KK^\dag K\Th
x,x\big\rangle-2\big\langle K\Th x,u\big\rangle
+\big\langle Ku,u\big\rangle\big)ds\Big]\\
\ns\ds=\frac{1}{2}\dbE\Big[\big\langle
P(t)\eta,\eta\big\rangle+\int_t^T\big(\big\langle
K\Th x,\Th x\big\rangle-2\big\langle K\Th x,u\big\rangle+\langle Ku,u\rangle\big)ds\Big]\\
\ns\ds=\frac{1}{2}\dbE\Big[\big\langle
P(t)\eta,\eta\big\rangle+\int_t^T \big\langle K(u-\Th x),u-\Th x\big\rangle ds\Big]\\
\ns\ds=\cJ\big(t,\eta;\Th(\cd)x(\cd)\big)
+\frac{1}{2}\dbE\int_t^T\big\langle K(u-\Th
x),u-\Th x\big\rangle ds.
\end{array}
$$

Hence,
$$\cJ(t,\eta;\Th(\cd)x(\cd))\leq \cJ(t,\eta;u),\q\forall u(\cd)\in\cU[t,T],$$
if and only if $P$ is a regular solution to the
Riccati equation {\rm(\ref{Riccati})}.

\ms

Proof of assertion ii). Without loss of
generality, we assume that $t=0$. The proof is
divided  into five steps.

\vspace{0.1cm}

{\bf Step 1}. In this step, we introduce some
operators and their finite dimensional
approximation.

Let $\cl\Th(\cd)$ be an optimal feedback
operator of Problem (SLQ) over $[0,T]$. For each
$s\in [0,T]$, define three operators $X_{s}$,
$Y_{s}$ and $\wt X_{s}$ on $H$ as follows:
\begin{equation}\label{8.20-eq1}
X_{s} \zeta =  \bar x(s;\zeta), \q Y_{s}\zeta =
\bar y(s;\zeta), \q \wt X_{s}\zeta = \tilde
x(s;\zeta), \q\forall\,\zeta\in H.
\end{equation}
For a.e. $s\in [0,T]$, define an operator
$Z_{s}$ on $H$ by
\begin{equation}\label{8.20-eq2}
Z_{s}\zeta = \bar z(s;\zeta), \qq
\forall\,\zeta\in H.
\end{equation}
Here $(\bar x,\bar y,\bar z)$ solves
\eqref{FBSDE5.1} and $\tilde x$ solves
\eqref{5.26-eq4}. Now we are going to show some
properties of the above four operators.

Denote by $I_n$ the identity matrix on $\dbR^n$
(or, the identity map on $H_n$). Consider the
following equations:
\begin{equation}\label{7.19-eq1}
\left\{
\begin{array}{ll}
\ds dX_n=  (A_n+A_{1,n}+B_n\Th_n)X_n ds+
(C_n+D_n\Th_n) X_n
dW(s)& \mbox{ in }[0,T],\\
\ns\ds dY_n=-\big[(A_n+A_{1,n})^* Y_n+ C^*_n
Z_n+ Q_nX_n\big]dt+ Z_ndW(t) &
\mbox{ in }[0,T],\\
\ns\ds X_n(0)=I_{n}, \q Y_n(T)= G_nX_n(T)
\end{array}
\right.
\end{equation}
and
\begin{equation}\label{7.19-eq2}
\left\{\2n
\begin{array}{ll}
\ds d\wt X_n\!=\!
\big[\!-A_n\!-\!A_{1,n}\!-\!B_n\Th_n\!+(C_n+D_n\Th_n)\big(C_n\!+D_n\Th_n\big)^\top
\big]^\top\wt X_n ds\\
\ns\ds\hspace{1.4cm}  - (C_n+D_n\Th_n)^\top\wt
X_n
dW(s)\hspace{4.4cm}\mbox{in }[0,T],\\
\ns\ds \wt X_n(0)=I_n.
\end{array}
\right.
\end{equation}
Clearly, both \eqref{7.19-eq1} and
\eqref{7.19-eq2} can be viewed as $\dbR^{n\times
n}\equiv\dbR^{n^2}$-valued equations. They admit
unique solutions
$ (X_n,Y_n,Z_n)\in C_\dbF([0,T];L^2(\Om;$
$\dbR^{n\times n}))\times C_\dbF([0,T];L^2(\Om;$
$\dbR^{n\times n}))\times
L^2_\dbF(0,T;\dbR^{n\times n}) $
and
$ \wt X_n \in C_\dbF([0,T];L^2(\Om;\dbR^{n\times
n}))$,
respectively. By It\^o's formula, we see that
$$
\begin{array}{ll}\ds
\q X_n(s)\wt X_n(s)^\top - X_n(0)\wt X_n(0)^\top\\
\ns\ds =\int_0^s
\big\{\langle(A_n\!+A_{1,n}\!+B_n\Th_n)X_n,\wt
X^\top\rangle + \langle X_n,\wt
X_n^\top\big[\!-A_n\!-A_{1,n}\!-B_n\Th_n\!+\big(C_n\!+\!D_n\Th_n\big)^2
\big]\rangle \\
\ns\ds\qq + \langle(C_n+D_n\Th_n) X_n, - \wt
X_n^\top(C_n+D_n\Th_n)\rangle dr\\
\ns\ds\q + \int_0^s\big[\langle (C_n+D_n\Th_n)
X_n, \wt X_n^\top\rangle + \langle X_n, - \wt
X_n^\top(C_n+D_n\Th_n)\rangle \big]dW(r)=0.
\end{array}
$$
Consequently,
\begin{equation}\label{8.20-eq5}
\wt X(s)_n^\top=X_n(s)^{-1},\qq
\dbP\mbox{-a.s.,}\q \forall\ s\in[0,T].
\end{equation}

Let $\zeta\in H$. It is an easy matter to see
that $x_n(s)=X_n(s)\G_n\zeta$,
$y_n(s)=Y_n(s)\G_n\zeta$,
$z_n(s)=Z_n(s)\G_n\zeta$ and $\tilde x_n(s)=\wt
X_n(s)\G_n\zeta$. Thus,
$(X_n(\cd)\G_n\zeta,Y_n(\cd)\G_n\zeta,Z_n(\cd)\G_n\zeta)$
solves \eqref{6.7-eq1} and $\wt
X_n(\cd)\G_n\zeta$ solves \eqref{6.7-eq2}.  For
each $s\in [0,T]$, define three operators
$X_{n,s}$, $Y_{n,s}$ and $\wt X_{n,s}$ on $H$ as
follows:
\begin{equation}\label{7.19-eq5}
X_{n,s}\G_n\zeta =  X_{n}(s)\G_n\zeta, \q
Y_{n,s}\G_n\zeta =  Y_{n}(s)\G_n\zeta, \q \wt
X_{n,s}\G_n\zeta = \wt X_{n}(s)\G_n\zeta,
\q\forall\,\zeta\in H.
\end{equation}
For a.e. $s\in [0,T]$, define an operator
$Z_{n,s}$ on $H_n$ by
\begin{equation}\label{7.19-eq5zz}
Z_{n,s}\G_n\zeta =  Z_{n}(s)\G_n\zeta, \qq
\forall\,\zeta\in H.
\end{equation}

In view of \eqref{7.19-eq5}, \eqref{7.19-eq5zz}
and \eqref{6.7-eq3}, we find that
\begin{equation}\label{7.24-eq4}
\left\{
\begin{array}{ll}\ds
\lim_{n\to+\infty}X_{n,s}\G_{n}\zeta =
X_s\zeta  &\mbox{ strongly in }L^2_{\cF_s}(\Om;H),\\
\ns\ds \lim_{n\to+\infty}Y_{n,s}\G_{n}\zeta =
Y_s\zeta  &\mbox{ strongly in }L^2_{\cF_s}(\Om;H),\\
\ns\ds\lim_{n\to+\infty}Z_{n,s}\G_{n}\zeta =
Z_s\zeta &\mbox{ strongly in }L^2_{\dbF}(0,T;H),\\
\ns\ds \lim_{n\to+\infty}\wt X_{n,s}\G_{n}\zeta
= \wt X_s \zeta &\mbox{ strongly in
}L^2_{\cF_s}(\Om;H).
\end{array}
\right.
\end{equation}

\medskip

{\bf Step 2}. In this step, we give an explicit
formula of $P(\cd)$.

By the well-posedness results for the equations
\eqref{6.7-eq1} and \eqref{6.7-eq2}, and the
fact that both $A$ and $-A^*$ generate
$C_0$-semigroups on $H$ (because $A$ generates a
$C_0$-group on $H$), we see that
$$
\begin{array}{ll}\ds
|X_{n,s}\G_n\zeta|_{L^2_{\cF_s}(\Om;H)}\leq
\cC|\zeta|_H,\q
|Y_{n,s}\G_n\zeta|_{L^2_{\cF_s}(\Om;H)}\leq
\cC|\zeta|_H,\\
\ns\ds
|Z_{n,\cd}\G_n\zeta|_{L^2_{\dbF}(0,T;H)}\leq
\cC|\zeta|_H,\q |\wt
X_{n,s}\G_n\zeta|_{L^2_{\cF_t}(\Om;H)}\leq
\cC|\zeta|_H,
\end{array}
$$
where the constant $\cC$ is independent of $n$.
This implies that
\begin{equation}\label{7.2016-eq2}
\begin{array}{ll}\ds
|X_{n,s}\G_n|_{\cL(H;L^2_{\cF_s}(\Om;H))}\leq
\cC,\q
|Y_{n,s}\G_n|_{\cL(H;L^2_{\cF_s}(\Om;H))}\leq
\cC,\\
\ns\ds
|Z_{n,\cd}\G_n|_{\cL(H;L^2_{\dbF}(0,T;H))}\leq
\cC,\q |\wt
X_{n,s}\G_n|_{\cL(H;L^2_{\cF_s}(\Om;H))}\leq
\cC.
\end{array}
\end{equation}

By \eqref{7.2016-eq2} and using \cite[Theorems
5.2-5.3]{LZ1}, we deduce that, there exist
subsequences $\{X_{n_k,s}\}_{k=1}^\infty\subset
\{X_{n,s}\}_{n=1}^\infty$,
$\{Y_{n_k,s}\}_{k=1}^\infty$ $\subset
\{Y_{n,s}\}_{n=1}^\infty$,
$\{Z_{n_k,s}\}_{k=1}^\infty\subset
\{Z_{n,s}\}_{n=1}^\infty$ and $\{\wt
X_{n_k,s}\}_{k=1}^\infty\subset \{\wt
X_{n,s}\}_{n=1}^\infty$ (these sequences may
depend on $s$), and (pointwise defined)
operators $X_1(s,\cd)$, $Y_1(s,\cd), \wt
X_1(s,\cd)\in \cL(H;L^2_{\cF_s}(\Om;H))$ (for
each $s\in [0,T]$) and $Z_1(\cd,\cd)\in
\cL(H;L^2_{\dbF}(0,T;$ $H))$  such that
\begin{equation}\label{7.24-eq3}
\left\{
\begin{array}{ll}\ds
\lim_{k\to+\infty}X_{n_k,s}\G_{n_k}\zeta =
X_1(s,\cd)\zeta &\mbox{ weakly in }L^2_{\cF_s}(\Om;H),\\
\ns\ds \lim_{k\to+\infty}Y_{n_k,s}\G_{n_k}\zeta
=
Y_1(t,\cd)\zeta &\mbox{ weakly in }L^2_{\cF_s}(\Om;H),\\
\ns\ds\lim_{k\to+\infty}Z_{n_k,s}\G_{n_k}\zeta =
Z_1(\cd,\cd)\zeta &\mbox{ weakly in }L^2_{\dbF}(0,T;H),\\
\ns\ds \lim_{k\to+\infty}\wt
X_{n_k,s}\G_{n_k}\zeta = \wt X_1(s,\cd)\zeta
&\mbox{ weakly in }L^2_{\cF_s}(\Om;H),
\end{array}
\right.
\end{equation}
and that
\begin{equation}\label{7.24-eq6}
\begin{array}{ll}\ds
|X_1(s,\cd)\zeta|_{L^2_{\cF_s}(\Om;H)}\leq
\cC|\zeta|_H,\q
|Y_1(s,\cd)\zeta|_{L^2_{\cF_s}(\Om;H)}\leq
\cC|\zeta|_H,\\
\ns\ds
|Z_1(\cd,\cd)\zeta|_{L^2_{\dbF}(0,T;H)}\leq
\cC|\zeta|_H,\q |\wt
X_1(s,\cd)\zeta|_{L^2_{\cF_s}(\Om;H)}\leq
\cC|\zeta|_H.
\end{array}
\end{equation}

It follows from \eqref{7.24-eq4} and
\eqref{7.24-eq3} that
\begin{equation}\label{8.20-eq3}
\begin{array}{ll}\ds
X(s,\cd)\=X_s(\cd) = X_1(s,\cd),\q
Y(s,\cd)\=Y_s(\cd) = Y_1(s,\cd),\\
\ns\ds Z(s,\cd)\=Z_s(\cd) = Z_1(s,\cd), \q \wt
X(t,\cd) \=\wt X_s(\cd)=\wt X_1(s,\cd).
\end{array}
\end{equation}

By \eqref{5.7-eq3}, \eqref{8.20-eq1},
\eqref{8.20-eq2} and \eqref{8.20-eq3}, we find
that
\begin{equation}\label{5.7-eq5.1}
R\Th X+B^* Y+D^* Z =0, \q \ae
(s,\om)\in[0,T]\times\Om.
\end{equation}

From \eqref{8.20-eq5}, \eqref{7.24-eq4},
\eqref{7.24-eq3} and \eqref{8.20-eq3}, it is
easy to see  that for any $s\in [0,T]$, $X(s)\wt
X(s)^*=I$, $\dbP$-a.s., that is, for any $s\in
[0,T]$, $\wt X(s)^*=X(s)^{-1}$, $\dbP$-a.s.

Put
\begin{equation}\label{8.20-eq9}
P(\cd)= Y(\cd) X(\cd)^{-1},\qq\Pi(\cd)= Z(\cd)
X(\cd)^{-1}.
\end{equation}
It follows from (\ref{5.7-eq5.1}) that
\begin{equation}\label{5.7-eq6.1}
B^* P +D^* \Pi +R \Th =0, \q \ae
(s,\om)\in[0,T]\times\Om.
\end{equation}

\medskip

{\bf Step 3}. In this step, we give an estimate
of the norm of $P(\cd)$ introduced in
\textbf{Step 2}.

 Let $s\in [0,T)$ and $\eta\in
L^2_{\cF_s}(\Om;H)$. Consider the following
FBSEE:
\begin{equation}\label{5.14-eq12}
\left\{
\begin{array}{ll}\ds
d x^s(r)=\big(A+A_1+B\Th\big) x^sdr + \big(C +
D \Th \big)x^sdW(r) &\mbox{\rm in } (s,T],\\
dy^s(r) = -\big[(A+A_1)^* y^s + C^* z^s + Qx^s
\big] dr +  z^sdW(r) &\mbox{\rm in } [s,T),
\\ \ns\ds x^s(s)=\eta, \q
y^s(T)= G x^s(T).
\end{array}
\right.
\end{equation}
By Lemma \ref{5.7-prop1}, it is easy to see that
\eqref{5.14-eq12} admits a unique solution
$
\big(x^s(\cd),y^s(\cd),z^s(\cd)\big)\big(\equiv\big(x^s(\cd;\eta),y^s(\cd;\eta),z^s(\cd;\eta)\big)\big)\in
C_\dbF([s,T];L^2(\Om;H))\times
C_\dbF([s,T];L^2(\Om;H))\times L^2_\dbF(s,T;H)$
such that
\begin{equation}\label{5.72016-eq3}
R\Th x^s(r)+B^* y^s(r)+D^* z^s(r) =0,\q\ae
(r,\om)\in (s,T)\times\Om.
\end{equation}
For every $r\in [s,T]$, define two families of
operators $X^s_r$ and $Y^s_r$ on
$L^2_{\cF_s}(\Om;H)$:
$$
X^s_r \eta \=  x^s(r;\eta), \q Y^s_r \eta  \=
y^s(r;\eta).
$$

It follows from Lemmas \ref{lm2} and \ref{lm3}
that for all $r\in [s,T]$ and $\eta\in
{L^2_{\cF_s}(\Om;H)}$,
\begin{equation}\label{5.26-eq8.1}
|X^s_r\eta|_{L^2_{\cF_r}(\Om;H)}\leq
\cC|\eta|_{L^2_{\cF_s}(\Om;H)},\q
|Y^s_r\eta|_{L^2_{\cF_r}(\Om;H)}\leq
\cC|\eta|_{L^2_{\cF_s}(\Om;H)}.
\end{equation}
This indicates that $X^s_r$ and $Y^s_r$ belong
to $\cL(L^2_{\cF_s}(\Om;H);L^2_{\cF_r}(\Om;H))$
for every $r\in [s,T]$.

By \eqref{5.14-eq12}, it is easy to see that,
for any $\zeta\in H$,
$$
X^s_rX(s)\zeta = x^s(r;X(s)\zeta)= x(r;\zeta).
$$
Thus,
$$
Y^s_s X(s)\zeta = y^s(s;X(s)\zeta)=Y(s)\zeta.
$$
This implies that
\begin{equation}\label{5.14-eq13}
Y^s_s =Y(s)\wt X(s)^*\q \mbox{ for all } s\in
[0,T],\;\ \dbP\hb{-a.s.}
\end{equation}

Let $\eta, \xi\in L^2_{\cF_s}(\Om;H)$. Since
$Y^s_r \eta=y^s(r;\eta)$ and $X^s_r
\xi=x^s(r;\xi)$, applying It\^o's formula to
$\langle y^s(\cd;\eta),x^s(\cd;\xi)\rangle$ and
noting \eqref{5.14-eq12}--\eqref{5.72016-eq3},
we obtain that
$$
\mE\langle GX^s_T \eta,X^s_T\xi \rangle  -
\mE\langle Y^s_s\eta, \xi \rangle  =
-\mE\int_s^T \big(\langle QX^s_r\eta,X^s_r\xi
\rangle  + \langle R\Th X^s_r\eta,\Th X^s_r\xi
\rangle \big)dr.
$$
Therefore,
$$
\begin{array}{ll}\ds
\mE\langle Y^s_s\eta, \xi \rangle  =
\mE\Big\langle (X^{s}_T)^* GX^s_T\eta
+\mE\int_s^T \big((X^{s}_r)^* Q X^s_r\eta  +
(X^{s}_r)^* \Th ^* R \Th X^s_r\eta \big)dr,\xi
\Big\rangle.
\end{array}
$$
We conclude from this that, for any $\eta\in
L^2_{\cF_s}(\Om;H)$,
\begin{equation}\label{5.14-eq11}
Y^s_s\eta  = \mE\((X^{s}_T)^* GX^s_T\eta
+\mE\int_s^T \big((X^{s}_r )^*Q X^s_r\eta  +
(X^{s}_r)^*\Th^* R \Th X^s_r\eta
\big)dr\;\Big|\;\cF_s\).
\end{equation}
It deduces that $Y(s)\wt X(s)^*=Y_s^s$ is
symmetric for any $s\in [0,T]$, $\dbP$-a.s.
Further, \eqref{5.14-eq11} together with
\eqref{5.26-eq8.1} implies that for any $s\in
[0,T]$ and $\eta \in L^2_{\cF_s}(\Om;H)$,
\begin{equation}\label{5.14-eq15}
\mE|Y^s_s\eta|_H^2 \leq \cC\mE|\eta|_H^2,
\end{equation}
where $\cC$ is  independent of $s\in [0,T)$.
According to \eqref{5.14-eq15}, we find that
\begin{equation}\label{5.14-eq16.1}
|Y(s)\wt
X(s)^*|_{\cL(L^2_{\cF_s}(\Om;H);\;L^2_{\cF_s}(\Om;H))}\leq
\cC.
\end{equation}
It follows from \eqref{8.20-eq9},
\eqref{5.14-eq13} and \eqref{5.14-eq16.1} that,
for some positive constant $\cC$,
\begin{equation}\label{5.14-eq9}
|P(s)|_{\cL(L^2_{\cF_s}(\Om;H);\;L^2_{\cF_s}(\Om;H))}\le
\cC, \qq\forall\;s\in [0,T].
\end{equation}

\medskip

{\bf Step 4}. In this step, we show that
$P(\cd)$ is a mild solution to an
operator-valued differential equation.

 Let $s\in [0,T)$ and $\zeta\in
L^2_{\cF_s}(\Om;H)$. Consider the following
FBSEE:
\begin{equation}\label{8.20-eq10}
\left\{
\begin{array}{ll}\ds
d x^s_n(\tau)=\big(A_n+A_{1,n}+B_n\Th_n\big)
x^s_nd\tau + \big(C_n +
D_n \Th_n \big)x_n^sdW(\tau) &\mbox{\rm in } [s, T],\\
dy^s_n(\tau) = -\big[(A_n+A_{1,n})^* y^s_n +
C_n^* z^s_n + Q_nx_n^s \big] d\tau +
z_n^sdW(\tau) &\mbox{\rm in } [s, T],
\\ \ns\ds x_n^s(s)=\zeta, \q
y_n^s(T)= G_n x_n^s(T).
\end{array}
\right.
\end{equation}
For every $r\in [s,T]$, define two families of
operators $X^{n,s}_r$ and $Y^{n,s}_r$ on
$L^2_{\cF_s}(\Om;H)$ as follows:
$$
X^{n,s}_r \zeta \=  x^s_n(r;\zeta), \q Y^{n,s}_r
\zeta \= y_n^s(r;\zeta).
$$
Similar to the proofs of \eqref{7.24-eq4} and
\eqref{5.14-eq13}, we can show that
\begin{equation}\label{8.20-eq11}
\left\{
\begin{array}{ll}\ds
\lim_{n\to+\infty}X^{n,s}_r\zeta =
X^{s}_r\zeta  &\mbox{ in }L^2_{\cF_r}(\Om;H),\\
\ns\ds \lim_{n\to+\infty}Y^{n,s}_r\zeta =
Y^{s}_r\zeta  &\mbox{ in }L^2_{\cF_r}(\Om;H)
\end{array}
\right.
\end{equation}
and
\begin{equation}\label{8.20-eq12}
Y^{n,s}_s =Y_n(s)\wt X_n(s)^\top\q \mbox{ for
all } s\in [0,T],\;\ \dbP\hb{-a.s.}
\end{equation}

Let $P_n(\cd)=Y_n(\cd)\wt X_n(\cd)^\top$ and
$\Pi_n(\cd)= Z_n(\cd) \wt X_n(\cd)^\top$.  We
conclude from \eqref{5.14-eq13},
\eqref{8.20-eq11} and \eqref{8.20-eq12} that
\begin{equation}\label{8.20-eq13}
\lim_{n\to+\infty}P_n(s)\G_n\eta = P(s)\eta
\q\mbox{ in } H,\qq \forall\, \eta\in H.
\end{equation}

By It\^o's formula,
$$
\begin{array}{ll}
\ds dP_n\3n&\ds=\Big\{\!-\big[(A_n+A_{1,n})^\top
Y_n\!+C_n^\top Z_n\!+ Q^\top_n X^\top_n\big] \wt
X^\top_n \!+ \! Y_n \wt X^\top_n\big[(C_n
+D_n\Th_n)^2\!-\!A_n\!-B_n\Th_n\big]\\
\ns&\ds\q
- Z_n \wt X^\top_n(C_n \!+D_n\Th_n)\Big\}d\tau +\big[ Z_n \wt X^\top_n- Y_n \wt X^\top_n(C_n+D_n\Th_n)\big]dW(\tau)\\
\ns&\ds =\Big\{\!-(A_n+A_{1,n})^\top
P_n\!-\!C_n^\top\Pi_n-Q_n +P_n\big[(C_n
+D_n\Th_n)^2-(A_n+A_{1,n})-B_n\Th_n\big]
\\
\ns&\ds\q-\Pi_n(C_n \!+D_n\Th_n)\Big\}d\tau
+\big[\Pi_n-P_n(C_n+D_n\Th_n)\big]dW(\tau).
\end{array}
$$
Let $\L_n=\Pi_n-P_n(C_n+D_n\Th_n)$.
Then
\begin{equation}\label{5.7-eq7}
\3n\begin{array}{ll} \ds
dP_n\3n&\ds=\!\Big\{-(A_n+A_{1,n})^\top
P_n-C_n^\top [\L_n+P_n(C_n+D_n\Th_n)]-Q_n
+P_n\big[(C_n\1n+D_n\Th)^2\\
\ns&\ds\q-(A_n+A_{1,n})-B_n\Th_n\big]
-[\L_n+P_n(C_n+D_n\Th_n)](C_n
+D_n\Th_n)\Big\}d\tau+\L_n
dW(\tau)\\
\ns&\ds = \big[ -  P_n(A_n+A_{1,n}) -
(A_n+A_{1,n})^\top P_n - \L_n C_n -
C_n^\top \L_n - C_n^\top P_nC_n\\
\ns&\ds\q - (P_nB_n +  C_n^\top P_nD_n + \L_n
D_n)\Th_n - Q_n\big]d\tau  + \L_n dW(\tau),
\end{array}
\end{equation}
and $P_n(T)=G_n$. This implies that
$(P_n(\cd),\L_n(\cd))$ is the adapted solution
to \eqref{5.7-eq7} with deterministic
coefficients and final datum. Thus, $P_n(\cd)$
is deterministic and $\L_n(\cd)=0$. Then,
\eqref{5.7-eq7} becomes
\begin{equation}\label{5.7-eq8}
\dot P_n+P_n(A_n+A_{1,n})+(A_n+A_{1,n})^\top
P_n+C_n^\top P_nC_n+(P_nB_n+C^\top_n P_nD_n
)\Th_n+Q_n=0.
\end{equation}
Thus,
$$
\begin{array}{ll}\ds
P_n(s)=e^{A_n^\top(T-s)}G_ne^{A_n(T-s)} +
\int_s^T e^{A_n^\top(r-s)}\big[P_nA_{1,n}
+A_{1,n}^\top
P_n+C_n^\top P_nC_n \\
\ns\ds\hspace{7.4cm}+ (P_nB_n + C^\top_n P_nD_n
)\Th_n + Q_n\big]e^{A_n(r-s)}dr.
\end{array}
$$
Therefore, for any $\zeta\in H$,
\begin{equation}\label{8.20-eq6}
\begin{array}{ll}\ds
P_n(s)\G_n\zeta\3n&\ds=e^{A_n^\top(T-s)}G_ne^{A_n(T-s)}\G_n\eta
+ \int_s^T e^{A_n^\top(r-s)}\big[P_nA_{1,n}
+A_{1,n}^\top P_n+C_n^\top
P_nC_n\\
\ns&\ds\hspace{5.14cm}+(P_nB_n\!+\!C^\top_n
P_nD_n )\Th_n\!+\!Q_n\big]e^{A_n(r-s)}\G_n\zeta
dr.
\end{array}
\end{equation}
From \eqref{8.20-eq8} and \eqref{8.20-eq13}, we
know that for any $\zeta\in H$,
\begin{equation}\label{8.20-eq14}
\lim_{n\to\infty}
e^{A_n^\top(T-s)}G_ne^{A_n(T-s)}\G_n\zeta
=e^{A^*(T-s)}Ge^{A(T-s)} \zeta
\end{equation}
and
\begin{equation}\label{8.20-eq15.1}
\begin{array}{ll}\ds
\lim_{n\to\infty}e^{A_n^\top(r-s)}\big[P_nA_{1,n}
+A_{1,n}^\top P_n+C_n^\top
P_nC_n+(P_nB_n+C^\top_n P_nD_n
)\Th_n+Q_n\big]e^{A_n(r-s)}\G_n\zeta\\
\ns\ds = e^{A^*(r-s)}\big[P A_{1} +A_{1}^*
P+C^\top PC+(PB+C^\top PD
)\Th+Q\big]e^{A(r-s)}\zeta.
\end{array}
\end{equation}
By \eqref{8.20-eq8} again, we have that
$$
\begin{array}{ll}\ds
\Big|e^{A_n^\top(r-s)}\big[P_nA_{1,n}
+A_{1,n}^\top P_n+C_n^\top
P_nC_n+(P_nB_n+C^\top_n P_nD_n
)\Th_n+Q_n\big]e^{A_n(r-s)}\G_n\zeta\Big|_H\\
\ns\ds\leq
\cC\big|e^{A_n^\top(r-s)}\big|_{\cL(H)}\big[2\big|A_{1,n}^\top\big|_{\cL(H)}
\big|P_n\big|_{\cL(H)}+\big|C_n^\top\big|_{\cL(H)}
\big|P_n\big|_{\cL(H)}\big|C_n\big|_{\cL(H)}+\big(\big|P_n\big|_{\cL(H)}\big|B_n\big|_{\cL(U;H)}\\
\ns\ds\qq+\big|C^\top_n\big|_{\cL(H)}
\big|P_n\big|_{\cL(H)}\big|D_n\big|_{\cL(U;H)}
)\big|\Th_n\big|_{\cL(H;U)}+\big|Q_n\big|_{\cL(H)}\big]\big|e^{A_n(r-s)}\big|_{\cL(H)}\big|\G_n\zeta\big|_H\\
\ns\ds \leq
\cC\big|e^{A^*(r-s)}\big|_{\cL(H)}\big[2\big|A_1^*\big|_{\cL(H)}
\big|P\big|_{\cL(H)}+\big|C^*\big|_{\cL(H)}
\big|P\big|_{\cL(H)}\big|C\big|_{\cL(H)}+\big(\big|P\big|_{\cL(H)}\big|B\big|_{\cL(U;H)}\\
\ns\ds\qq+\big|C^*\big|_{\cL(H)}
\big|P\big|_{\cL(H)}\big|D\big|_{\cL(U;H)}
)\big|\Th\big|_{\cL(H;U)}+\big|Q\big|_{\cL(H)}\big]\big|e^{A(r-s)}\big|_{\cL(H)}\big|\zeta\big|_H.
\end{array}
$$
This, together with \eqref{8.20-eq15.1} and
Lebesgue's dominated convergence theorem,
implies that
\begin{equation}\label{8.20-eq15}
\begin{array}{ll}\ds
\lim_{n\to\infty}\int_s^T\!e^{A_n^\top(r-s)}\big[P_nA_{1,n}
\!+\!A_{1,n}^\top P_n\!+\!C_n^\top
P_nC_n\!+(P_nB_n\!+\!C^\top_n P_nD_n
)\Th_n\!+\!Q_n\big]e^{A_n(r-s)}\G_n\zeta ds\\
\ns\ds = \int_s^Te^{A^*(r-s)}\big[PA_{1}
+A_{1}^* P_n+C^* PC+(PB+C^\top PD
)\Th+Q\big]e^{A(r-s)}\zeta dr.
\end{array}
\end{equation}
It follows from \eqref{8.20-eq13},
\eqref{8.20-eq6}, \eqref{8.20-eq14} and
\eqref{8.20-eq15} that
\begin{equation}\label{8.20-eq16}
\begin{array}{ll}\ds
P(s)\zeta =e^{A^*(T-s)}Ge^{A(T-s)}\zeta +
\int_s^T
e^{A^*(r-s)}\big[PA_{1} + A_{1}^* P+C^* PC\\
\ns\ds\hspace{6.95cm}+(PB+ C^* PD )\Th +
Q\big]e^{A(r-s)}\zeta dr.
\end{array}
\end{equation}

{\bf Step 5}. Finally, in this step, we prove
that $P(\cd)$ solves the Riccati equation
\eqref{Riccati}.

 From (\ref{5.7-eq6.1}), we see
that
\begin{equation}\label{5.7-eq9}
0=B^* P+D^* P(C+D\Th)+R\Th =B^* P+D^* PC
+K\Th,\q\ae (s,\om)\in [0,T]\times\Om.
\end{equation}
Consequently,
\begin{equation}\label{5.7-eq9.1.1}
0= \Th^*B^* P+\Th^*D^* PC +\Th^*K\Th,\q\ae
(s,\om)\in [0,T]\times\Om.
\end{equation}
Using (\ref{5.7-eq9.1.1}), \eqref{8.20-eq16} can
be written as
$$
\begin{array}{ll}\ds
P(t)\zeta=e^{(T-s)A}Ge^{(T-s)A}\zeta + \int_s^T
e^{(r-s)A}\big[PA_{1} \!+\!A_{1}^* P+P B\Th \\
\ns\ds\hspace{4.14cm}+ \Th^*B^*P+(C+D\Th)^*
P(C+D\Th)+\Th^* R\Th+Q\big]e^{(r-s)A}\zeta ds.
\end{array}
$$
Since $P(T)=G\in\dbS(H)$, and $Q(\cd)$ and
$R(\cd)$ are symmetric operator-valued
functions, by Proposition \ref{prop3}, we have
$P(\cd)\in C_\cS([0,T];\dbS(H))$. Thus $K(\cd)$
is a symmetric operator-valued function. Hence,
for a.e. $s\in [0,T]$, $K(s)$ admits a
generalized pseudo inverse $K(s)^\dag$.

Since $\L_n=0$ for all $n\in\dbN$, we see that
$$
\Pi_n-P_n(C_n+D_n\Th_n)=Z_n\wt
X_n^{\top}-P_n(C_n+D_n\Th_n)=0 \q\mbox{ for
every }n\in\dbN.
$$
Therefore,
$$
\Pi-P(C+D\Th)=0.
$$
This implies that
$$
(B^* P+D^* PC)=- K\Th.
$$
Thus, \eqref{regular-1} holds and
$$
K^\dag(B^* P+D^* PC)=-K^\dag K\Th.
$$
Noting that $K^\dag K$ is an orthogonal
projection, we see that \eqref{regular-2} holds
and
$$
\Th=-K^\dag(B^* P+D^* PC)+\big(I-K^\dag K\big)\th
$$
for some $\th(\cd)\in L^2(0,T;\cL(H;U))$.
Consequently,
\begin{equation}\label{5.7-eq9.1}
\ds(PB+C^* PD )\Th=\Th^* KK^\dag(B^* P+D^* PC
)=-(PB+C^* PD)K^\dag(B^* P+D^* PC).
\end{equation}
From \eqref{5.7-eq9.1} and \eqref{8.20-eq16}, we
obtain the Riccati equation \eqref{Riccati}.
This completes the proof of the ``only if" part.
\endpf

\ms


\section{Proof of Theorem \ref{th4.6}}\label{sec-pr-main2}


In this section, we prove Theorem \ref{th4.6}.
Without loss of generality, we assume that
$t=0$.

{\it Proof of Theorem \ref{th4.6} }: (i) $\Ra$
(ii). The proof is long. We divide it into three
steps.

{\bf Step 1}. In this step, we introduce a
sequence of operator-valued functions
$\{P_j\}_{j=1}^N$.

Let $P_0$ be the solution to
\begin{equation}\label{8.20-eq30}
\left\{\2n
\begin{array}{ll}
\ds\dot P_0+P_0(A+A_1)+(A+A_1)^* P_0+C^* P_0C+Q=0 &\mbox{ in }[0,T),\\
\ns\ds P_0(T)=G.
\end{array}
\right.
\end{equation}
Applying Proposition \ref{prop4.5} to
\eqref{8.20-eq30} with $\Th=0$, we obtain that
\begin{equation}\label{8.20-eq30.1}
R(s)+D(s)^* P_0(s)D(s)\geq\l I,\q P_0(s)\geq\a
I,\qq\ae~s\in[0,T].
\end{equation}
Inductively, for $j = 0,1,2, \cdots$, we set
\begin{equation}\label{Iteration-i}
\begin{array}{ll}\ds
K_j\=R+D^* P_jD,\qq L_j\=B^* P_j+D^* P_jC,
\\
\ns\ds \Th_j\=-K_j^{-1}L_j,
\qq\cA_j\=A_1+B\Th_j,\qq C_j\=C+D\Th_j,
\end{array}
\end{equation}
and let $P_{j+1}$ be the solution to
\begin{equation}\label{}
\left\{\2n
\begin{array}{ll}
\ds\dot{P}_{j+1}+P_{j+1}(A+\cA_j)+(A+\cA_j)^*
P_{j+1}+C_j^* P_{j+1}C_j
+\Th_j^* R\Th_j+Q=0 &\mbox{ in }[0,T),\\
\ns\ds P_{j+1}(T)=G.
\end{array}
\right.
\end{equation}

{\bf Step 2}. In this step, we show the uniform
boundedness of the sequence
$\{P_j\}_{j=1}^\infty$.

From \eqref{8.20-eq30.1}, we have that
\begin{equation}\label{11.27-eq1}
K_0(s)\geq\l I,\q P_0(s)\geq\a
I,\qq\ae~s\in[0,T].
\end{equation}
Then $\Th_0=-K_0^{-1}L_0\in L^2(0,T;\cL(H;U))$.
It follows from Proposition \ref{prop4.5} (with
$P$ and $\Th$ in \eqref{eq-Lya} replaced by
$P_{1}$ and $\Th_0$, respectively) that
$$
K_{1}(s)\geq\l I,\q P_{1}(s)\geq\a I,\qq
\ae~s\in[0,T].
$$
Inductively, we have that
\begin{equation}\label{R+Pi-lowerbound}
K_{j+1}(s)\geq\l I,\q P_{j+1}(s)\geq\a I,\qq
\ae~s\in[0,T],\qq j=0,1,2,\cdots.
\end{equation}
We now claim that $\{P_j\}_{j=1}^\i$ converges
uniformly in $C_\cS([0,T];\dbS(H))$. To show
this, let
$$\D_j\deq P_j-P_{j+1},\qq \Upsilon_j\deq\Th_{j-1}-\Th_j,\qq j\geq1.$$
Then for $j\geq1$ and $\zeta\in H$, we have
\begin{equation}\label{Di-equa1}
\begin{array}{ll}\ds
\q -\D_j(s)\zeta=P_{j+1}(s)\zeta-P_j(s)\zeta \\
\ns\ds= \int_s^T
e^{(r-s)A^*}\big[\big(P_j(s)-P_{j+1}(s)\big)
\cA_j +
\cA_j^* \big(P_j(s)-P_{j+1}(s)\big)+ C^*_j \big(P_j(s)-P_{j+1}(s))C_j\\
\ns\ds\qq\qq\qq +P_j(\cA_{j-1}-\cA_j)+
(\cA_{j-1}-\cA_j)^*P_j+C_{j-1}^*
P_jC_{j-1}-C_j^* P_jC_j+\Th_{j-1}^*
R\Th_{j-1}\\
\ns\ds\qq\qq\qq-\Th_j^* R\Th_j\big]e^{(r-s)A}\zeta ds\\
\ns\ds= \int_s^T e^{(r-s)A^*}\big[\D_j(s) \cA_j
+ \cA_j^* \D_j(s)+ C^*_j
\D_j(s)C_j+P_j(\cA_{j-1}-\cA_j)+
(\cA_{j-1}-\cA_j)^*P_j\\
\ns\ds\qq\qq\qq +C_{j-1}^* P_jC_{j-1}-C_j^*
P_jC_j+\Th_{j-1}^* R\Th_{j-1} -\Th_j^*
R\Th_j\big]e^{(r-s)A}\zeta dr
\end{array}
\end{equation}
From \eqref{Iteration-i}, we have that
$$
\cA_{j-1}-\cA_j= A_1+B\Th_{j-1} - A_1- B\Th_j =
B(\Th_{j-1} - \Th_j) =B\Upsilon_j,
$$
$$
C_{j-1}-C_j = C+D\Th_{j-1} - C -
D\Th_j=D(\Th_{j-1} -  \Th_j) =D\Upsilon_j,
$$
and
$$
\begin{array}{ll}\ds
C_{j-1}^* P_jC_{j-1}-C_j^* P_jC_j \\
\ns\ds= (C+D\Th_{j-1})^*P_j(C+D\Th_{j-1}) -
(C+D\Th_{j})^*P_j(C+D\Th_{j})\\
\ns\ds = \Th_{j-1}^*D^*\!P_jD\Th_{j-1}\!\! -\!
\Th_{j}^*D^*P_jD\Th_{j}\! +\!
\Th_{j-1}^*D^*\!P_jC \!+\!
C^*\!P_jD\Th_{j-1}\!\!-\!\Th_{j}^*\!D^*\!P_jC
\!-\!
C^*P_jD\Th_{j}\\
\ns\ds =
(\Th_{j-1}-\Th_{j})^*D^*P_jD(\Th_{j-1}-\Th_{j})+
(C+D\Th_{j})^* P_jD(\Th_{j-1}-\Th_{j})\\
\ns\ds\q +(\Th_{j-1}-\Th_{j})^* D^*
P_j(C+D\Th_{j})
\\ \ns\ds =\Upsilon_j^* D^* P_jD\Upsilon_j+C_j^*
P_jD\Upsilon_j+\Upsilon_j^* D^* P_jC_j.
\end{array}
$$
Similarly, we can obtain that
\begin{equation*}\label{Di-equa2}
\left\{\2n
\begin{array}{ll}
\ds  \Th_{j-1}^* R\Th_{j-1}-\Th_j^*
R\Th_j=\Upsilon_j^* R\Upsilon_j+\Upsilon_j^*
R\Th_j+\Th_j^* R\Upsilon_j,\\
\ns\ds B^* P_j+D^* P_jC_j+R\Th_j=B^* P_j+D^*
P_jC+(R+D^* P_jD)\Th_j=0.
\end{array}
\right.
\end{equation*}
These, together with \eqref{Di-equa1}, yields
that
\begin{eqnarray}\label{Di-equa3}
&&  \D_j(s)-\int_s^Te^{(r-s)A^*}\big(\Delta_j\cA_j+\cA_j^*\D_j+C_j^*\D_jC_j\big)e^{(r-s)A}dr \nonumber\\
&&=\int_s^Te^{(r-s)A^*}\big(P_j B\Upsilon_j+\Upsilon_j^* B^* P_j+\Upsilon_j^* D^* P_jD\Upsilon_j+C_j^* P_jD\Upsilon_j+\Upsilon_j^* D^* P_jC_j\nonumber\\
&&\qq\qq\qq+\Upsilon_j^* R\Upsilon_j+\Upsilon_j^* R\Th_j+\Th_j^* R\Upsilon_j\big)e^{(r-s)A}ds\\
&&=\!\int_s^T\!\!e^{(r-s)A^*}\!\big[\Upsilon_j^*K_j\Upsilon_j
\!+\!(P_jB\!+\!C_j^* P_jD\!+\!\Th_j^*
R)\Upsilon_j\!+\!\Upsilon_j^*(B^* P_j\!+\!D^*
P_jC_j\!+\!R\Th_j)\big]e^{(r-s)A}dr\nonumber
\\
&&=\int_s^Te^{A^*(r-s)}\Upsilon_j^*K_j\Upsilon_je^{A(r-s)}dr.\nonumber
\end{eqnarray}
From \eqref{Di-equa3}, we know that $\D_j(\cd)$
is a solution to \eqref{lm2.4-eq1} with $\wt
G=0$, $\wt A = \cA_j$, $\wt C = \cC_j$ and $\wt
Q = \Upsilon_j^*K_j\Upsilon_j\geq 0$. Using
Lemma \ref{lm2.4}, we have that $\D_j(\cd)\geq
0$, namely, $P_{j-1}(\cd)-P_j(\cd)\geq 0$ for
$j\geq1$. Noting (\ref{R+Pi-lowerbound}), we
obtain
$$P_1(s)\geq P_j(s)\geq P_{j+1}(s)\geq\a I,\qq\forall s\in [0,T],\q\forall j\geq1.$$
Therefore, the sequence $\{P_j\}_{j=1}^\i$ is
uniformly bounded. Consequently, there exists a
constant $\cC>0$ such that (noting
(\ref{R+Pi-lowerbound})) for all $j\geq0$ and
a.e. $s\in [0,T]$,
\begin{equation}\label{Di-equa4}
\left\{\2n
\begin{array}{ll}
\ns\ds|P_j(s)|_{\cL(H)}\leq \cC,\qq
|K_j(s)|_{\cL(U)}\leq \cC,
\\
\ns\ds |\Th_j(s)|_{\cL(H;U)}\leq \cC\big(|B(s)|_{\cL(U;H)}+|C(s)|_{\cL(H)}\big),\\
\ns\ds|\cA_j(s)|_{\cL(H)}\leq
|A_1(s)|_{\cL(H)}+\cC|B(s)|_{\cL(U;H)}\big(|B(s)|_{\cL(U;H)}+|C(s)|_{\cL(H)}\big),\\
\ns\ds |C_j(s)|_{\cL(H)}\leq
\cC\big(|B(s)|_{\cL(U;H)}+|C(s)|_{\cL(H)}\big).
\end{array}
\right.
\end{equation}

{\bf Step 3}. In this step, we prove the
convergence of the sequence
$\{P_j\}_{j=1}^\infty$.

Noting
$$
\L_j=\Th_{j-1}-\Th_j
=K_j^{-1}D^*\D_{j-1}DK_{j-1}^{-1}L_j-K_{j-1}^{-1}\big(B^*\D_{j-1}+D^*\D_{j-1}C\big)
$$
and (\ref{Di-equa4}), one has
\begin{equation}\label{3.22}
\begin{array}{ll}
\ds|\Upsilon_j(s)^*
K_j(s)\Upsilon_j(s)|_{\cL(H)}\\
\ns\ds\leq\big(|\Th_j(s)|_{\cL(U)}+|\Th_{j-1}(s)|_{\cL(U)}\big)
|K_j(s)|_{\cL(U)}
|\Th_{j-1}(s)-\Th_j(s)|_{\cL(U)}\\
\ns\ds \leq
\cC\big(|B(s)|_{\cL(U;H)}+|C(s)|_{\cL(H)}\big)^2|\D_{j-1}(s)|_{\cL(H)}.
\end{array}
\end{equation}
Equation (\ref{Di-equa3})  implies that
$$
\D_j(s)=-\int^T_s e^{(r-s)A^*}\big(\D_j\cA_j+
\cA_j^*\D_j+C_j^*\D_iC_j+\Upsilon_j^*
K_j\Upsilon_j\big)e^{(r-s)A}dr.
$$
Making use of (\ref{3.22}) and noting
(\ref{Di-equa4}), we get
$$
|\D_j(s)|_{\cL(H)}\leq
\int^T_s\f(r)\big(|\D_j(r)|_{\cL(H)}
+|\D_{j-1}(r)|_{\cL(H)}\big)dr,\qq\forall
s\in[0,T],\q\forall j\geq1,$$
where $\f(\cd)$ is a nonnegative integrable
function independent of $\D_j(\cd)$. By
Gronwall's inequality,
$$
|\D_j(s)|_{\cL(H)} \leq
e^{\int_0^T\f(r)dr}\int^T_s\f(r)|\D_{j-1}(r)|_{\cL(H)}dr\equiv
b\int^T_s\f(r)|\D_{j-1}(r)|_{\cL(H)}dr,
$$
where $b=e^{\int_0^T\f(r)dr}$. Set
$a\deq\max_{0\leq r\leq T}|\D_0(r)|_{\cL(H)}$.
By induction, we deduce that
$$
|\D_j(s)|_{\cL(H)}\leq
a\frac{b^j}{j!}\(\int_s^T\f(r)dr\)^j,\qq\forall
s\in[0,T],
$$
which implies the uniform convergence of
$\{P_j\}_{j=1}^\i$. Denote by $P$ the limit of
$\{P_j\}_{j=1}^\infty$, then (noting
(\ref{R+Pi-lowerbound}))
$$
K(s)=\lim_{j\to\i}K_j(s) \geq\l I,
\qq\ae~s\in[0,T],$$
and as $j\to\infty$,
$$\left\{\2n\begin{array}{ll}
\ds \Th_j\to-K^{-1}L
\equiv\Th\q\hb{in } L^2(0,T;\cL(H;U)),\\
\ns\ds \cA_j\to A_1+B\Th\q\hb{in }
L^1(0,T;\cL(H)),\qq C_j\to C+D\Th\q\hb{in
}L^2(0,T;\cL(H)).\end{array}\right.$$
Therefore, $P(\cd)$ solves the following
equation (in the sense of mild solution):
$$\left\{\2n
\begin{array}{ll}
\ds\dot P+P(A+B\Th)+(A+B\Th)^* P+(C+D\Th)^* P(C+D\Th)+\,\Th^* R\Th+Q=0 &\mbox{ in }[0,T),\\
\ns\ds P(T)=G,
\end{array}
\right.$$
which is equivalent to (\ref{Riccati}).

\ms

(ii) $\Ra$ (i). Let $P(\cd)$ be the strongly
regular solution to (\ref{Riccati}). Then there
exists a $\l>0$ such that
\begin{equation}\label{iitoi}
K(s)\geq \l I,\qq\ae~s\in[0,T].
\end{equation}
Set $\Th\deq-K^{-1}L\in L^2(0,T;\cL(H;U))$. For
any $u(\cd)\in\cU[0,T]$, let
$x(\cd)=x(\cd;0,0,u)$ be the solution to
\eqref{state} with $t=0$ and $\eta=0$. Applying
It\^o's formula to $s\mapsto\langle
P(s)x(s),x(s)\rangle$, we have
$$
\begin{array}{ll}
\ds \cJ(0,0;u(\cd))=\dbE\[\lan Gx(T),x(T)\ran
+\int_0^T\big(\lan Qx,x\ran + \lan Ru,u\ran\big)ds\]\\
\ns\ds=\dbE\int_0^T\Big\{\lan-\big[P\big(A+A_1\big)
+\big(A+A_1\big)^* P
+C^* PC+Q-L^* K^\dag L\big]x,x\ran\\
\ns\ds\qq\qq ~+\lan P\big[\big(A+A_1\big)x+Bu\big],x\ran+\lan Px,\big(A+A_1\big)x+Bu\ran\\
\ns\ds\qq\qq~\2n+\lan P\big(Cx+Du\big),Cx+Du\ran+\lan Qx,x\ran+\lan Ru,u\ran\Big\}ds  \\
\ns\ds=\dbE\int_0^T\big(\lan\Th^*K\Th x,x\ran
-2\lan K\Th x,u\ran +\lan Ku,u\ran\big)ds\\
\ns\ds=\dbE\int_0^T\lan K\big(u-\Th x\big),u-\Th
x\ran ds.
\end{array}$$
Noting (\ref{iitoi}) and making use of Lemma
\ref{lm2.5}, we obtain that for some $c_0>0$ and
all $u(\cd)\in\cU[0,T]$,
$$
\begin{array}{ll}\ds
\cJ(0,0;u(\cd)) \ds =\dbE\int_0^T\lan
K\big(u-\Th x\big),u-\Th x\ran ds \geq c_0\l
\dbE\int_0^T|u(s)|^2ds.
\end{array}
$$
Then (i) holds. \endpf

\ms

\begin{remark}\label{rmk4.7}
From the first part of the proof of Theorem
\ref{th4.6}, we see that if (\ref{J>l*}) holds,
then the strongly regular solution to
(\ref{Riccati}) satisfies (\ref{strong-regular})
with the same constant $\l>0$.
\end{remark}
%


\section{The uniform convexity of the cost functional}\label{sec-convex}


In this section, we study the uniform convexity
of the cost functional.

Define four operators as follows:
$$
\G_t:H\to \cX[t,T]\=L^2_{\dbF}(t,T;H),\q
\G_t\eta=x(\cd\,;t,\eta,0), \qq \forall\,
\eta\in H,
$$
where $x(\cd\,;t,\eta,0)$ is the solution to
\eqref{state} with $u\equiv0$;

$$
\Xi_t:\cU[t,T]\to \cX[t,T],\q \Xi_t
u=x(\cd\,;t,0,u), \qq \forall\, u\in \cU[t,T],
$$
where $x(\cd\,;t,0,u)$ is the solution to
\eqref{state} with $\eta=0$;

$$
\widehat\G_t:H\to L^2_{\cF_T}(\Om;H),\q \widehat
\G_t\eta=x(T\,;t,\eta,0), \qq \forall\, \eta\in
H;
$$
and
$$
\widehat \Xi_t:\cU[t,T]\to L^2_{\cF_T}(\Om;H),\q
\widehat \Xi_tu=x(T\,;t,0,u), \qq \forall\,
\eta\in H.
$$

From the well-posedness of \eqref{state}, we
find that  all these operators are bounded
linear operators. Then, for any $t\in[0,T)$ and
$(\eta,u(\cd))\in H\times\cU[t,T]$, the
corresponding state process $x(\cd)$ and its
terminal value $x(T)$ are given by
$$
x(\cd)=(\G_t\eta)(\cd)+(\Xi_tu)(\cd),\q
x(T)=\widehat\G_t\eta+\widehat \Xi_tu.
$$
Hence, the cost functional can be written as
\begin{equation}\label{5.18-eq1}
\begin{array}{ll}
\ds \cJ(t,\eta;u(\cd))=\lan
G\big(\widehat\G_t\eta+\widehat
\Xi_tu\big),\widehat\G_t\eta+\widehat \Xi_tu\ran
+\lan Q(\G_t\eta+\Xi_tu),\G_t\eta+\Xi_tu\ran
+\langle Ru,u\rangle.
\end{array}
\end{equation}

Recall that $u(\cd)\mapsto \cJ(t,\eta;u(\cd))$
is uniformly convex if and only if for some
$\l>0$,
\begin{equation}\label{J>l}
\cJ(t,0;u(\cd))\geq\l\dbE\int_t^T|u(s)|^2ds,\qq\forall
u(\cd)\in\cU[t,T].
\end{equation}
From \eqref{5.18-eq1},  \eqref{J>l} is
equivalent to the following:
\begin{equation}\label{M_2>l}
\widehat \Xi_t^*G\widehat
\Xi_t+\Xi_t^*Q\Xi_t+R\geq\l I,\q \mbox{ for some
$\l>0$.}
\end{equation}

According to the above argument, we get the
following result:
\begin{proposition}\label{th-covvex}
The map $u(\cd)\mapsto \cJ(t,\eta;u(\cd))$ is
uniformly convex if and only if \eqref{M_2>l}
holds.
\end{proposition}

It is obvious that if the condition
\eqref{classical} hold, then \eqref{M_2>l} holds
for $\l=\d$. On the other hand, if $R \geq \d I$
does not hold, \eqref{M_2>l} still may be true
if  $G$ is large enough.  An example is given
below:

We first introduce the following assumption.

\vspace{0.1cm}

{\bf (AS3)} Assume that $D$ is invertible for
a.e. $t\in [0,T]$ and $D(\cd)^{-1}\in L^\infty
(0,T;\cL(H;U))$.

\vspace{0.1cm}

Without loss of generality, we set $t=0$ in
\eqref{state}. Under  {\bf (AS3)}, the control
system \eqref{state} can be written as a BSEE
\begin{equation}\label{state1}
\left\{\2n
\begin{array}{ll}
\ns\ds dx =\big[(A+A_1 )x +B D^{-1}\hat x
\big]ds+\big(Cx+\hat x\big)dW(s) &\mbox{ in }[0,T), \\
\ns\ds x(T) \in L^2_{\cF_T}(\Om;H),
\end{array}
\right.
\end{equation}
where $\hat x(\cd)=D(\cd)u(\cd)$ and $x(T)$ is
the value of the solution to \eqref{state} at
time $T$. If {\bf (AS1)} and {\bf (AS3)} hold,
then \eqref{state1} is well-posed and
$$
|(x(\cd),\hat
x(\cd))|_{C_\dbF([0,T];L^2(\Om;H))\times
L^2_\dbF(0,T;H)}\leq
\cC|x(T)|_{L^2_{\cF_T}(\Om;H)}.
$$
This implies that there is a $\cC_0>0$,
depending only on $A$, $A_1$, $B$, $C$ and $D$
such that
\begin{equation}\label{8.20-eq31}
|u(\cd)|^2_{L^2_\dbF(0,T;H)}\leq
\cC_0|x(T)|^2_{L^2_{\cF_T}(\Om;H)}.
\end{equation}

Let us make the following further assumption.

\vspace{0.1cm}

{\bf (AS4)} There is a
$\mu_0>\cC_0\big(|R|_{L^\infty(0,T;\cL(U))}+\e_0\big)$
with $\e_0>0$ such that for any $\zeta\in H$,
$\lan G\zeta,\zeta\ran\geq \mu_0|\zeta|^2$.

\vspace{0.1cm}

If {\bf (AS4)} holds,  it follows from
\eqref{8.20-eq31} that
$$
\lan\widehat \Xi_t^*G\widehat \Xi_t u,u\ran =
\lan G\widehat \Xi_t u,\h \Xi_t u\ran \geq
\mu_0|u|_{L^2_\dbF(0,T;H)}^2 \geq
\big(|R|_{L^\infty(0,T;\cL(U))}+\e_0\big)|u|_{L^2_\dbF(0,T;H)}^2.
$$
This deduces that \eqref{M_2>l} holds for
$\l=\e_0$.

According to the above argument, we have the
following result.

\begin{theorem}
Let {\bf (AS1)}, {\bf (AS3)} and {\bf (AS4)}
 hold. Then the map $u(\cd)\mapsto
\cJ(t,\eta;u(\cd))$ is uniformly convex.
\end{theorem}
%


\section{Appendix: Proofs for some preliminary results}


\subsection{Proof of Lemma \ref{lm2}}

{\it Proof of Lemma \ref{lm2}}\,:  Without loss
of generality, we assume that $t=0$. Write
$$
N=\lceil\frac{1}{\e}\big(\big|\cA\big|^2_{L^1(0,T;\cL(H))}
+ \big|\cB\big|^2_{L^2(0,T;\cL(H))}\big)\rceil
+1,
$$
where $\e>0$ is a constant to be determined
later. Define a sequence of
$\{\tau_{j,\e}\}_{j=1}^N$ as follows:
$$
\left\{
\begin{array}{ll} \ds
\tau_{1,\e} = \Big\{r\in
[0,T]\;\Big|\;\(\int_0^r|\cA(s)|_{\cL(H)}ds\)^2+\int_0^r|\cB(s)|^2_{\cL(H)}ds=\e\Big\},\\
\ns\ds \tau_{k,\e} =  \Big\{r\in
[0,T]\;\Big|\;\(\int_{\tau_{k-1,\e}}^r|\cA(s)|_{\cL(H)}ds\)^{2}+\int_{\tau_{k-1,\e}}^r|\cB(s)|^2_{\cL(H)}ds=\e\Big\},\\
\ns\ds\qq\qq\; k= 2,\cdots,N.
\end{array}
\right.
$$

Consider the following  SEE:
\begin{equation}\label{10.29-eq9}
\left\{
\begin{array}{ll}\ds
dx = \big(A x + \tilde f\big)ds + \tilde gdW(s) &\mbox{ in }(0,T],\\
\ns\ds x(0)=\eta.
\end{array}
\right.
\end{equation}
Here $\tilde f\in L^2_\dbF(\Om;L^1(0,T;H))$ and
$\tilde g\in L^2_\dbF(0,T;H)$. Clearly,
\eqref{10.29-eq9} admits a unique mild solution
$x\in C_\dbF([0,T];L^2(\Om;H))$ (e.g.
\cite[Chapter 6]{Prato}). Define a map
$\cJ:\;C_\dbF([0,\tau_{1,\e}];L^2(\Om;H))\to
C_\dbF([0,\tau_{1,\e}];L^2(\Om;H))$ as follows:
$$
C_\dbF([0,\tau_{1,\e}];L^2(\Om;H))\ni\tilde
x\mapsto x=\cJ(\tilde x),
$$
where $x$ is the solution to \eqref{10.29-eq9}
with $\tilde f$ and $\tilde g$ replaced by $\cA
\tilde x + f$ and $\cB \tilde x + g$,
respectively. We claim that $\cJ$ is
contractive. Indeed, for any $\tilde x_1,\tilde
x_2\in C_\dbF([0,\tau_{1,\e}];L^2(\Om;H))$,
\begin{equation}\label{10.29-eq11}
\begin{array}{ll}
\q|\cJ(\tilde x_1) - \cJ(\tilde
x_1)|_{C_\dbF([0,\tau_{1,\e}];L^2(\Om;H))}^2  \\
\ns\ds \leq 2\!\sup_{s\in
[0,\tau_{1,\e}]}\mE\(\Big|\int_0^s e^{(s-r)A}\cA
(\tilde x_1-\tilde x_2)dr\Big|_H\)^2 + 2\mE
\Big| \int_0^{\tau_{1,\e}} \big|e^{(s-r)A}\cB
(\tilde
x_1-\tilde x_2)\big|_H^2dW(r)\Big|_H \\
\ns\ds \leq 2\(\sup_{s\in
[0,T]}|e^{As}|^2_{\cL(H)}\)\big||\cA|_{\cL(H)}\big|^2_{L^1(0,\tau_{1,\e})}|\tilde
x_1-\tilde x_2|^2_{C_\dbF([0,\tau_{1,\e}];L^2(\Om;H))}\\
\ns\ds\q + 2\(\sup_{s\in
[0,T]}|e^{As}|^2_{\cL(H)}\)\big||\cB|_{\cL(H)}\big|^2_{L^2(0,\tau_{1,\e})}|\tilde
x_1-\tilde
x_2|^2_{C_\dbF([0,\tau_{1,\e}];L^2(\Om;H))}.
\end{array}
\end{equation}
Let
\begin{equation}\label{2.9-eq4}
M_T=\sup_{s\in[0,T]}|e^{As}|_{\cL(H)},\qq \e =
\frac{1}{16M_T^2}.
\end{equation}
It follows from \eqref{10.29-eq11} that
\begin{equation}\label{10.29-eq12}
\big|\cJ(\tilde x_1) - \cJ(\tilde
x_1)\big|_{C_\dbF([0,\tau_{1,\e}];L^2(\Om;H))}^2
\leq \frac{1}{4}\big|\tilde x_1-\tilde
x_2\big|^2_{C_\dbF([0,\tau_{1,\e}];L^2(\Om;H))}.
\end{equation}
This shows that $\cJ$ is contractive. Hence, it
has a unique fixed point $x\in
C_\dbF([0,\tau_{1,\e}];L^2(\Om;H))$, which
solves \eqref{6.20-eq1} in $[0,\tau_{1,\e}]$ (in
the sense of mild solution). Inductively, we
conclude that \eqref{6.20-eq1} admits a mild
solution $x$ in $[\tau_{k-1,\e},\tau_{k,\e}]$
for $k=2,\cds,N$. Furthermore,
\begin{eqnarray}\label{10.29-eq13}
&&\ds\q|x|_{C_\dbF([0,\tau_{1,\e}];L^2(\Om;H))}^2 \nonumber \\
&&\ds \leq 4\sup_{s\in
[0,\tau_{1,\e}]}\mE\(\Big|\int_0^s e^{A(t-r)}\cA
x dr\Big|_H\)^2 + 4\sup_{t\in
[0,\tau_{1,\e}]}\mE\(\Big| \int_0^s
e^{A(s-r)}\cB
x dW(r)\Big|_H\)^2 \\
&&\ds \q + 4 \sup_{s\in [0,T]}\mE
\(\Big|e^{As}\eta  + \int_0^s e^{A(s-r)}f dr +
\int_0^t e^{A(s-r)}g dW (r)\Big|_H\)^2\nonumber\\
&&\ds \leq
4M_T^2\big||\cA|_{\cL(H)}\big|^2_{L^1(0,\tau_{1,\e})}|
x|^2_{C_\dbF([0,\tau_{1,\e}];L^2(\Om;H))} +
4M_T^2\big||\cB|_{\cL(H)}\big|^2_{L^2(0,\tau_{1,\e})}|
x|^2_{C_\dbF([0,\tau_{1,\e}];L^2(\Om;H))}\nonumber\\
&&\ds \q + \cC\big(|\eta|_{H}^2 +
|f|_{L^2_\dbF(\Om;L^1(0,T;H))}^2 +
|g|_{L^2_\dbF(0,T;H)}^2\big).\nonumber
\end{eqnarray}
This, together with the choice of $\tau_{1,\e}$,
implies that
\begin{equation}\label{10.29-eq14}
|x|_{C_\dbF([0,\tau_{1,\e}];L^2(\Om;H))}^2\leq
\cC\big(|\eta|_{H}^2 +
|f|_{L^2_\dbF(\Om;L^1(0,T;H))}^2 +
|g|_{L^2_\dbF(0,T;H)}^2\big).
\end{equation}
Repeating the above argument, we obtain
\eqref{lm2-eq1}. The uniqueness of the solution
is obvious.
\endpf


\subsection{Proof of Lemma \ref{lm3}}


{\it Proof of Lemma \ref{lm3}}\,: Without loss
of generality, we assume that $t=0$.

Write \vspace{-0.3cm}
$$
N=\lceil\frac{1}{\e}\big(\big|A_1^*\big|^2_{L^1(t,T;\cL(H))}
+ \big|C^*\big|^2_{L^2(t,T;\cL(H))}\big)\rceil
+1,
$$
where $\e>0$ is a constant to be determined
later. Define a sequence of
$\{\tau_{k,\e}\}_{k=1}^N$ as follows:
$$
\left\{\3n
\begin{array}{ll} \ds
\tau_{1,\e} = \Big\{r\in
[0,T]\;\Big| \(\int_r^T|A_1(s)^*|_{\cL(H)}ds\)^2+\int_r^T|C(s)^*|^2_{\cL(H)}ds=\e\Big\},\\
\ns\ds \tau_{k,\e} =  \Big\{r\in [0,T]\;\Big|
\(\int_r^{\tau_{k-1,\e}}\!\!|A_1(s)^*|_{\cL(H)}ds\)^{2}+\!\int_r^{\tau_{k-1,\e}}\!\!|C(s)^*|^2_{\cL(H)}ds=\e\Big\},
\end{array} k= 2,\cdots,N.
\right.
$$

Consider the following  BSEE:
\begin{equation}\label{2.9-eq1}
\left\{
\begin{array}{ll}\ds
dy = -(A^*y+\tilde h) ds + zdW(s) &\mbox{ in }[0,T),\\
\ns\ds y(T)=\xi,
\end{array}
\right.
\end{equation}
where $\tilde h\in L^2_\dbF(\Om;L^1(0,T;H))$.
Clearly, \eqref{2.9-eq1} admits a unique mild
solution $(y(\cd),z(\cd))\in
C_\dbF([0,T];L^2(\Om;H))\times L^2_\dbF(0,T;H)$
(e.g. \cite{Mahmudov1}) such that
\begin{equation}\label{12092717-eq4}
\begin{array}{ll}\ds
y(s) = e^{A^*(T-s)} \xi + \int_s^T
e^{A^*(r-s)}\tilde h(r)dr - \int_s^T
e^{A^*(r-s)}z(r)dW(r),
\end{array}
\end{equation}
and
\begin{equation}\label{12.17-eq5}
z(s)= e^{A^*(T-s)}l(s)+\int_s^T e^{A^*(\eta-s)}
\k(\eta,s)d\eta,
\end{equation}
where $l(\cd)\in L^2_\dbF(0,T;H)$ such that
\begin{equation}\label{12.17-eq1}
\mE(\xi\;|\;\cF_s)=\mE \xi + \int_0^s
l(\si)dW(\si).
\end{equation}
and $\k(\cd,\cd)\in L^1(0,T;L_\dbF^2(0,T;H))$
such that for a.e. $s\in [0,T]$,
\begin{equation}\label{12.17-eq2}
\tilde h(s) =\mE \tilde h(s) + \int_0^s
\k(s,\si)dW(\si)
\end{equation}
and
\begin{equation}\label{12.17116-eq2}
|\k(\cd,\cd)|_{L^1(r,T;L_\dbF^2(r,T;H))}\le
|\tilde h(\cd)|_{L^1_\dbF(r,T;L^2(\Omega;H))}.
\end{equation}

Define a map
$\cJ:\;L^2_\dbF(\Om;C([\tau_{1,\e},T];H))\times
L^2_\dbF(\tau_{1,\e},T;H)\to
L^2_\dbF(\Om;C([\tau_{1,\e},T];H))\times
L^2_\dbF(\tau_{1,\e},T;H)$ as follows:
$$
L^2_\dbF(\Om;C([\tau_{1,\e},T];H))\times
L^2_\dbF(\tau_{1,\e},T;H)\ni(\tilde y,\tilde
z)\mapsto (y,z)=\cJ(\tilde y,\tilde z),
$$
where $(y,z)$ is the solution to \eqref{2.9-eq1}
with $\tilde h$ replaced by $A_1^* \tilde y +
C^*\tilde z+h$. We claim that $\cJ$ is
contractive. Indeed, for $j=1,2$ and $(\tilde
y_j,\tilde z_j) \in
L^2_\dbF(\Om;C([\tau_{1,\e},T];H))\times
L^2_\dbF(\tau_{1,\e},T;H)$, denote by $\tilde
h_j=A_1^* \tilde y_j + C^*\tilde z_j+h$ and
$\k_j(\cd,\cd)\in L^1(0,T;L_\dbF^2(0,T;H))$ such
that for a.e. $s\in [0,T]$,
\begin{equation}\label{12.17-eq2}
\tilde h_j(s) =\mE \tilde h_j(s) + \int_0^s
\k_j(s,\si)dW(\si)
\end{equation}
and
\begin{equation}\label{12.17116-eq2}
|\k_j(\cd,\cd)|_{L^1(r,T;L_\dbF^2(r,T;H))}\le
|\tilde h_j(\cd)|_{L^1_\dbF(r,T;L^2(\Omega;H))}.
\end{equation}
Let $\wt
M_T=\sup_{s\in[0,T]}|e^{A^*s}|_{\cL(H)}$. Then
\begin{equation}\label{2.9-eq2}
\begin{array}{ll}\ds
\q|\tilde z_1- \tilde
z_2|_{L^2_\dbF(\tau_{1,\e},T;H)}^2\\
\ns\ds = \mE\int_{\tau_{1,\e}}^T\Big|\int_s^T
e^{A^*(r-s)}
\big(\k_1(r,s)-\k_2(r,s)\big)dr\Big|^2_H
ds\\
\ns\ds  \leq \wt M_T^2\[\(\int_{\tau_{1,\e}}^T
|A_1(s)^*|_{\cL(H)}  ds\)^2 +
\int_{\tau_{1,\e}}^T |C(s)^*|_{\cL(H)}^2
ds\]\\
\ns\ds\q\times\big(\big|\tilde y_1 -\tilde y_2
\big|_{L^2_\dbF(\Om;C([\tau_{1,\e},T];H))}^2
+\big|\tilde z_1 -\tilde z_2
\big|_{L^2_\dbF(\tau_{1,\e},T;H)}^2 \big)\\
\ns\ds\leq \wt M_T^2 \e\big(\big|\tilde y_1
-\tilde y_2
\big|_{L^2_\dbF(\Om;C([\tau_{1,\e},T];H))}^2
+\big|\tilde z_1 -\tilde z_2
\big|_{L^2_\dbF(\tau_{1,\e},T;H)}^2 \big)
\end{array}
\end{equation}
and
\begin{eqnarray}\label{2.9-eq3}
&&|\tilde y_1- \tilde
y_2|_{L^2_\dbF(\Om;C([\tau_{1,\e},T];H))}^2\\
&& \leq 2\mE\Big|\int_{\tau_{1,\e}}^T\big|
e^{A^*(s-\tau_{1,\e})}\big(\tilde h_1(s)-\tilde
h_2(s)\big)\big|_Hds \Big|^2 + 2\mE
\int_{\tau_{1,\e}}^T
\big|e^{A^*(s-\tau_{1,\e})}\big(z_1(s)-z_2(s)\big)\big|_H^2ds\\
&&\leq 2\wt M_T^2|\tilde h_1(\cd)-\tilde
h_2(\cd)|^2_{L^1_\dbF(\tau_{1,\e},T;L^2(\Omega;H))}
+ 2\wt M_T^2|\tilde z_1- \tilde
z_2|_{L^2_\dbF(\tau_{1,\e},T;H)}^2\\
&&\leq 2(\wt M_T^2+\wt M_T^4)\e\big(\big|\tilde
y_1 -\tilde y_2
\big|_{L^2_\dbF(\Om;C([\tau_{1,\e},T];H))}^2
+\big|\tilde z_1 -\tilde z_2
\big|_{L^2_\dbF(\tau_{1,\e},T;H)}^2 \big).
\end{eqnarray}
Thus,
\begin{equation}\label{2.9-eq5}
\begin{array}{ll}
\q|\cJ(\tilde y_1,\tilde z_1) - \cJ(\tilde y_2, \tilde z_2)|_{L^2_\dbF(\Om;C([\tau_{1,\e},T];H))\times L^2_\dbF(\tau_{1,\e},T;H)}^2\\
\ns\ds \leq (3\wt M_T^2+2\wt
M_T^4)\e\big(\big|\tilde y_1 -\tilde y_2
\big|_{L^2_\dbF(\Om;C([\tau_{1,\e},T];H))}^2
+\big|\tilde z_1 -\tilde z_2
\big|_{L^2_\dbF(\tau_{1,\e},T;H)}^2 \big).
\end{array}
\end{equation}
Let  $\e = \frac{1}{4(3\wt M_T^2+2\wt M_T^4)}$.
It follows from \eqref{2.9-eq5} that
\begin{equation}\label{2.9-eq6}
\begin{array}{ll}\ds
|\cJ(\tilde y_1,\tilde z_1) - \cJ(\tilde y_2,
\tilde
z_2)|_{L^2_\dbF(\Om;C([\tau_{1,\e},T];H))\times
L^2_\dbF(\tau_{1,\e},T;H)}^2\\
\ns\ds \leq \frac{1}{4}\big|(\tilde y_1,\tilde
z_1)- (\tilde y_2, \tilde
z_2)\big|_{L^2_\dbF(\Om;C([\tau_{1,\e},T];H))\times
L^2_\dbF(\tau_{1,\e},T;H)}^2.
\end{array}
\end{equation}
Consequently, $\cJ$ is contractive. Hence, it
has a unique fixed point $(y,z)\in
L^2_\dbF(\Om;C([\tau_{1,\e},T];$ $H))\times
L^2_\dbF(\tau_{1,\e},T;H)$, which solves
\eqref{6.20-eq10} in $[\tau_{1,\e},T]$ (in the
sense of mild solution). Inductively, we
conclude that \eqref{6.20-eq10} admits a mild
solution $(y,z)$ in
$[\tau_{k,\e},\tau_{k-1,\e}]$ for $k=2,\cds,N$.
Furthermore,
\begin{equation}\label{2.9-eq7}
\begin{array}{ll}
\q|z|_{L^2_\dbF(\tau_{1,\e},T;H)}^2\\
\ns\ds \leq \mE\int_{\tau_{1,\e}}^T
\Big|e^{A^*(T-s)}l(s)+\int_s^T e^{A^*(r-s)}
\k(r,s)dr\Big|_H^2 ds\\
\ns\ds \leq 4\mE\int_{\tau_{1,\e}}^T
\big|e^{A^*(T-s)}l(s)\big|_H^2 ds +
\frac{3}{2}\mE\int_{\tau_{1,\e}}^T \Big|\int_s^T
e^{A^*(r-s)} \k(r,s)dr\Big|_H^2 ds\\
\ns\ds \leq \cC|\xi|_{L^2_{\cF_T}(\Om;H)}^2 +
\frac{3}{2}\wt M_T^2
\Big\{\int_{\tau_{1,\e}}^T\big[\mE \big|A_1(s)^*
y(s) + C(s)^*z(s) + h(s)
\big|_H^2\big]^{\frac{1}{2}}ds\Big\}^2\\
\ns\ds \leq \cC|\xi|_{L^2_{\cF_T}(\Om;H)}^2 +
\cC|h|_{L^2_{\dbF}(\Om;L^1(\tau_{1,\e},T;H))}^2
+ 2\wt M_T^2\e\big(|y
|_{L^2_\dbF(\Om;C([\tau_{1,\e},T];H))}^2 +|z
|_{L^2_\dbF(\tau_{1,\e},T;H)}^2 \big)
\end{array}
\end{equation}
and
\begin{eqnarray}\label{2.9-eq8}
&&\q|y|_{L^2_\dbF(\Om;C([\tau_{1,\e},T];H))}^2\nonumber\\
&&\ds \leq \mE\[\sup_{t\in
[\tau_{1,\e},T]}\Big|e^{A^*(T-t)} \xi - \int_t^T
e^{A^*(s-t)}\big(A_1^* \tilde y + C^*\tilde
z+h\big)ds - \int_t^T
e^{A^*(s-t)}z(s)dW(s)\Big|_H\]^2\nonumber\\
&&\ds \leq \cC|\xi|_{L^2_{\cF_T}(\Om;H)}^2 +
2\mE\(\sup_{t\in [\tau_{1,\e},T]}\Big|\int_t^T
e^{A^*(s-t)}\big(A_1^* \tilde y + C^*\tilde
z\big) dr\Big|_H\)^2 + 2 \mE
\int_{\tau_{1,\e}}^T \big|e^{A^*(s-t)}
z\big|_H^2 ds \nonumber\\
&&\ds \q + \cC|h|_{L^2_{\dbF}(\Om;L^1(0,T;H))}^2 \\
&&\ds \leq \cC|\xi|_{L^2_{\cF_T}(\Om;H)}^2 +
\cC|h|_{L^1_{\dbF}(0,T;L^2(\Om;H))}^2+ 4\wt
M_T^2|y|_{L^2_\dbF(\Om;C([\tau_{1,\e},T];H))}^2+
6\wt
M_T^2|z|_{L^2_\dbF(\tau_{1,\e},T;H)}^2.\nonumber
\end{eqnarray}
It follows from \eqref{2.9-eq7}, \eqref{2.9-eq8}
and the choice of $\tau_{1,\e}$  that
\begin{equation}\label{2.9-eq9}
|(y,z)|_{L^2_\dbF(\Om;C([\tau_{1,\e},T];H))\times
L^2_\dbF(\tau_{1,\e},T;H)}^2\leq
\cC\big(|\xi|_{L^2_{\cF_T}(\Om;H)}^2 +
|h|_{L^1_\dbF(0,T;L^2(\Om;H))}^2 \big).
\end{equation}
Repeating the above argument, we obtain
\eqref{2.9-eq0}. The uniqueness of the solution
is obvious.
\endpf


\subsection{Proof of Lemma \ref{5.7-prop1}}


{\it Proof of Lemma \ref{5.7-prop1}}\,: By
Lemmas \ref{lm2} and \ref{lm3}, the equation
\eqref{FBSDE5.1} admits a mild solution $(\bar
x(\cd), \bar y(\cd),\bar z(\cd))$. Since
$\cl\Th(\cd)$ is an optimal feedback operator of
Problem (SLQ), $(\bar x,\Th \bar x)$ is an
optimal pair. Put $\bar u(\cd)=\Th \bar x(\cd)$.
Fix arbitrarily a control $u(\cdot)\in
L^2_\dbF(t,T;U)$ and put
$$
u^\e(\cdot) = \bar u(\cdot) + \e [u(\cdot) -
\bar u(\cdot)] = (1-\e)\bar u(\cdot) + \e
u(\cdot) \in L^2_\dbF(t,T;U), \q\forall\;\e \in
[0,1].
$$
Denote by $x^\e(\cdot)$ the solution to
\eqref{state} corresponding to the control
$u^\e(\cdot)$. Write $\ds x_1^\e(\cd) =
\frac{1}{\e}\big[x^\e(\cd)-\bar x(\cd)\big]$ and
$\d u(\cd) = u(\cd) - \bar u(\cd)$. It is easy
to see that $x_1^\e(\cdot)$ solves the following
SEE:
\begin{equation}\label{fsystem3x}
\left\{
\begin{array}{lll}\ds
dx_1^\e = \big[(A+A_1)x_1^\e  + B \d u \big]ds +
\big(C x^\e_1 + D \d u \big)dW(s) &\mbox{ in
}(t,T],\\
\ns\ds x_1^\e(t)=0.
\end{array}
\right.
\end{equation}
Thanks to that $(\bar x(\cdot),\bar u(\cdot))$
is an optimal pair of Problem (SLQ), it holds
that
\begin{equation}\label{var 1}
\begin{array}{ll}\ds
0\leq \lim_{\e\to 0}\frac{\cJ(t,\eta;u^\e(\cdot)) - \cJ(t,\eta;\bar u(\cdot))}{\e} \\
\ns\ds\;\;\;= \dbE\int_t^T \big(\big\langle
Q\bar x,x_1^\e\big\rangle  + \big\langle R\bar
u,\d u\big\rangle\big) ds + \dbE\big\langle
G\bar x(T),x_1^\e(T)\big\rangle.
\end{array}
\end{equation}
By It\^o's formula,
\begin{equation}\label{max eq1}
\dbE \big\langle G\bar x(T),x_1^\e(T)\big\rangle
+  \dbE \int_t^T \big\langle Q\bar
x,x_1^\e\big\rangle ds= \dbE \int_t^T
\big(\big\langle B \d u, \bar y\big\rangle +
\big\langle D\d u, \bar z\big\rangle \big)ds.
\end{equation}
Combining \eqref{var 1} and \eqref{max eq1}, we
find that for any $u(\cdot)\in L^2_\dbF(0,T;U)$,
\begin{equation}\label{max ine2}
\dbE\int_t^T\big\langle R \bar u + B^*\bar y  +
D^*\bar z, u-\bar u \big\rangle ds\geq 0.
\end{equation}
This implies  \eqref{5.7-eq3}.
\endpf


\subsection{Proof of Lemma \ref{lm4}}


{\it Proof of Lemma \ref{lm4}}\,: Without loss
of generality, we assume that $t=0$.

In the beginning, we prove the first equality in
\eqref{6.7-eq3}. Let $\e\in [0,T]$ such that
\begin{equation}\label{2.9-eq15}
\e=\sup\Big\{ r\in [0,T]:10M_T^2\[\(\int_0^r
|A_1 |_{\cL(H)}ds \)^2 + \int_0^r |C|_{\cL(H)}^2
ds\]\leq \frac{1}{2}\Big\}\wedge T.
\end{equation}
For any $r\in [0,\e]$, it follows from
\eqref{FBSDE5.1} and \eqref{6.7-eq1} that
\begin{eqnarray}\label{2.9-eq14} &&\ds
\3n\3n\3n\q\sup_{r\in[0,\e]}\mE\big(|x(r)- x_n(r) |_H\big)^2 \nonumber \\
&&\ds\3n\3n\3n\leq
5\(\sup_{r\in[0,\e]}\mE\big|\big(e^{At} \!-\!
e^{At}\G_n\big)\zeta\big|_H^2\! +\!
\sup_{r\in[0,\e]}\mE\(\Big|\!\int_0^r\!
\big(e^{A(r-s)}A_1 x -\! e^{A(r-s)}\G_n A_{1,n}
x_n \big)ds\Big|_H\)^2  \\
&&\ds\qq \3n\3n\3n\!\!+\!
\sup_{r\in[0,\e]}\mE\(\!\Big|\!\int_0^r\!\!
\big(e^{A(r-s)} Bu\!\!-\!e^{A(r-s)}\G_n
B_nu\big) ds\Big|_H\)^2\!\! +\! \int_0^\e\!\!
\big|e^{A(r-s)}C x\! -\!
e^{A(r-s)}\G_n C_{n} x_n \big|_H^2ds\nonumber\\
&&\ds\qq\3n\3n\3n\!\!+ \int_0^\e \big|e^{A(r-s)}
Du-e^{A(r-s)}\G_n D_nu\big|_H^2ds\).\nonumber
\end{eqnarray}
Let us estimate the terms in the right hand side
of \eqref{2.9-eq14}. Clearly,
\begin{eqnarray}\label{2.9-eq10}
&&\ds
\3n\3n\3n\3n\3n\!\!\q\sup_{r\in[0,\e]}\mE\[\Big|
\int_0^r \big(e^{A(r-s)}A_1 x - e^{A(r-s)}\G_n
A_{1,n}
x_n \big)ds\Big|_H\]^2  \\
&&\ds\3n\3n \3n\3n\3n\!\!\leq \!\cC \mE
\Big(\!\int_0^\e |\big(A_1 - \G_n
A_{1}\!\G_n\big) x|_Hds \Big)^2   +\! 2M_T^2
\(\!\int_0^\e\! |A_{1}|_{\cL(H)} ds\)^2\!
\sup_{r\in[0,\e]}\mE\big|x(r) -\!
x_n(r)\big|_H^2,\nonumber
\end{eqnarray}
\begin{equation}\label{2.9-eq11}
\begin{array}{ll}\ds \sup_{r\in[0,\e]}\mE\Big| \int_0^r
\big[e^{A(r-s)} Bu\!-\!e^{A(r-s)}\G_n B_nu\big]
ds\Big|_H^2\leq \cC\(\int_0^\e
|(B-B_n)u|_Hds\)^2,
\end{array}
\end{equation}
\begin{equation}\label{2.9-eq12}
\begin{array}{ll}\ds \q \mE\int_0^\e
\big|e^{A(r-s)}C x - e^{A(r-s)}\G_n
C_{n} x_n \big|_H^2ds\\
\ns\ds \leq\cC\mE\int_0^\e\! \big|(C - C_n) x
\big|_H^2ds +\! 2M_T^2 \(\int_0^\e
|C|_{\cL(H)}^2 ds\)
\sup_{r\in[0,\e]}\mE\big|x(r) - x_n(r)\big|_H^2,
\end{array}
\end{equation}
and
\begin{equation}\label{2.9-eq13}
\begin{array}{ll}\ds \mE\int_0^\e
\big|e^{A(r-s)} Du-e^{A(r-s)}\G_n
D_nu\big|_H^2ds \leq \cC\int_0^\e \big|(D-
D_n)u\big|_H^2ds.
\end{array}
\end{equation}
Combining \eqref{2.9-eq14}, \eqref{2.9-eq11},
\eqref{2.9-eq12} and \eqref{2.9-eq13}, and
noting \eqref{2.9-eq15}, we get that
\begin{equation}\label{2.9-eq16}
\begin{array}{ll}\ds
\ds
\q\sup_{r\in[0,\e]}\mE|x(r)- x_n(r) |_H^2  \\
\ns\ds \leq \cC\mE\Big\{\big|\big(e^{Ar} -
e^{Ar}\G_n\big)\zeta\big|_H^2 + \Big[\int_0^\e
\big|\big(A_1 - \G_n
A_{1}\G_n\big) x \big|_Hds \Big]^2 \\
\ns\ds\qq\; + \[\int_0^\e |(B-B_n)u|_Hds\]^2
+\int_0^\e \big|(C-C_n) x \big|_H^2ds  +
\int_0^\e \big|(D- D_n)u\big|_H^2ds\Big\}.
\end{array}
\end{equation}
From \eqref{6.8-eq5.1}, we know that
\bel{zx-s1} \lim_{n\to\infty}\big|\big(e^{Ar} -
e^{Ar}\G_n\big)\zeta\big|_H^2=0. \ee
Since\vspace{-0.3cm}
$$
\Big[\int_0^\e \big|\big(A_1 - \G_n
A_{1}\G_n\big) x \big|_H ds \Big]^2 \leq
\[\int_0^\e|A_1|_{\cL(H)}|x|_Hds\]^2,
$$
we get from Lebesgue's dominated convergence
theorem that
\begin{equation}\label{2.9-eq17}
\begin{array}{ll}\ds
\lim_{n\to+\infty}\mE\Big[\int_0^\e
\big|\big(A_1 - \G_n A_{1}\G_n\big) x\big|_H ds
\Big]^2 =\mE\lim_{n\to+\infty}\Big[\int_0^\e
\big|\big(A_1 - \G_n A_{1}\G_n\big) x\big|_H ds
\Big]^2.
\end{array}
\end{equation}
Noting that \vspace{-0.2cm}
$$
\big|e^{A(r - s)}\big(A_1 - \G_n A_{1}\G_n\big)
x\big|_H \leq M_T|A_1|_{\cL(H)}|x|_H,
$$
by Lebesgue's dominated convergence theorem,
\eqref{2.9-eq17} and \eqref{6.8-eq5}, we have
that
\begin{equation}\label{2.9-eq18}
\begin{array}{ll}
\ds \lim_{n\to+\infty}\mE\Big[\int_0^\e
\big|\big(A_1 - \G_n A_{1}\G_n\big) x\big|_H ds
\Big]^2 =\mE\Big[\int_0^\e \lim_{n\to+\infty}
\big|\big(A_1 - \G_n A_{1}\G_n\big) x\big|_Hds
\Big]^2=0.
\end{array}
\end{equation}
By a similar argument, we can prove that
\begin{equation}\label{2.9-eq19}
\ds \lim_{n\to+\infty} \mE\Big\{\[\!\int_0^\e\!
|(B - B_n)u|_Hds\]^2\! +\!\int_0^\e\! \big|(C -
C_n) x \big|_H^2ds \!+ \!\int_0^\e \big|(D-
D_n)u\big|_H^2ds\Big\}\!=\!0.
\end{equation}

It follows from \eqref{2.9-eq16},
\eqref{2.9-eq18} and \eqref{2.9-eq19} that
$$
\begin{array}{ll}\ds
\lim_{n\to+\infty}\sup_{r\in[0,\e]}\mE|x(r)-
x_n(r) |_H^2=0.
\end{array}
$$
Repeating the above argument yields the first
equality in \eqref{6.7-eq3}. The proof for the
fourth one is similar. Now we consider the
second and third one.

Let
\begin{equation}\label{2.9-eq40}
T_1\=\inf\Big\{r\in [0,T]: \big[24(14\wt
M_T^2\!+\!1)\wt M_T^2 \!+\! 14\wt
M_T^2\big]\[\(\int_{r}^T \!\big|
A_{1}^*\big|_{\cL(H)}ds\)^2\! + \!\int_{r}^T\!
\big| C^*\big|_{\cL(H)}^2ds\] \leq
\frac{1}{2}\Big\}.
\end{equation}
We first recall that
\begin{equation}\label{2.9-eq20}
\begin{array}{ll}\ds
y(r) = e^{A^*(T-r)}G x(T) - \int_r^T
e^{A^*(s-r)}( A_{1} ^* y+ C^*z+ Qx)ds - \int_r^T
e^{A^*(s-r)}z(s)dW(s),
\end{array}
\end{equation}
and
\begin{equation}\label{2.9-eq21}
z(r)= e^{A^*(T-r)}l(r)+\int_r^T e^{A^*(\si-r)}
\k(\si,r)d\si,
\end{equation}
where $l(\cd)\in L^2_\dbF(0,T;H)$ such that
\begin{equation}\label{2.9-eq22}
\mE(G x(T)|\;\cF_r)=G\mE  x(T) + \int_0^r
l(\si)dW(\si)
\end{equation}
and $\k(\cd,\cd)\in L^1(0,T;L_\dbF^2(0,T;H))$
such that for a.e. $s\in [0,T]$,
\begin{equation}\label{2.9-eq23}
\big(A_{1} ^* y+ C^*z+ Qx\big)(s) =\mE
\big(A_{1} ^* y+ C^*z+ Qx\big)(s) + \int_0^s
\k(s,\si)dW(\si)
\end{equation}
and
\begin{equation}\label{2.9-eq24}
|\k(\cd,\cd)|_{L^1(r,T;L_\dbF^2(r,T;H))}\le
|A_{1} ^* y+ C^*z+
Qx|_{L^1_\dbF(r,T;L^2(\Omega;H))}.
\end{equation}
Furthermore,\vspace{-0.2cm}
\begin{equation}\label{2.9-eq25}
\begin{array}{ll}\ds
y_n(r) = e^{A_n^*(T-r)}G_n x_n(T) - \int_r^T
e^{A_n^*(s-r)}( A_{1,n} ^* y_n+ C_n^*z_n+ Q_nx_n)ds\\
\ns\ds\qq\qq - \int_r^T
e^{A_n^*(s-r)}z_n(s)dW(s),
\end{array}
\end{equation}
and\vspace{-0.2cm}
\begin{equation}\label{2.9-eq26}
z_n(r)= e^{A_n^*(T-r)}l_n(r)+\int_r^T
e^{A_n^*(\si-r)} \k_n(\si,r)d\si,
\end{equation}
where $l_n(\cd)\in L^2_\dbF(0,T;H)$ such that
\begin{equation}\label{2.9-eq27}
\mE(G_n x_n(T)|\;\cF_r)=G_n\mE  x_n(T) +
\int_0^r l_n(\si)dW(\si).
\end{equation}
and $\k_n(\cd,\cd)\in L^1(0,T;L_\dbF^2(0,T;H))$
such that for a.e. $s\in [0,T]$,
\begin{equation}\label{2.9-eq28}
\big(A_{1,n} ^* y_n+ C_n^*z_n+ Q_nx_n\big)(s)
=\mE \big(A_{1,n}^* y_n+ C_n^*z_n+
Q_nx_n\big)(s) + \int_0^s \k_n(s,\si)dW(\si)
\end{equation}
and
\begin{equation}\label{2.9-eq29}
|\k_n(\cd,\cd)|_{L^1(r,T;L_\dbF^2(r,T;H))}\le
|A_{1,n}^* y_n+ C_n^*z_n+
Q_nx_n|_{L^1_\dbF(r,T;L^2(\Omega;H))}.
\end{equation}
Then,
\begin{eqnarray}\label{2.9-eq30}
&&\ds \mE\sup_{r\in [T_1,T]}|y(r)-y_n(r)|_{H}^2
+ (14\wt
M_T^2+1)\mE\int_{T_1}^T|z(r)-z_n(r)|_{H}^2dr \nonumber\\
&&\leq \cC\mE\sup_{r\in
[T_1,T]}\[\big|\big(e^{A^*(T-r)}-e^{A_n^*(T-r)}\big)G
x(T)\big|_{H}^2+
 \big|e^{A_n^*(T-r)}(G-G_n)x(T)\big|_{H}^2
\nonumber\\
&&\ds\q +  \big|e^{A_n^*(T-r)}G_n
(x(T)-x_n(T))\big|_{H}^2 +  \(\!\int_{r}^T\!
\big|\big(e^{A^*(s-r)}\!-\!e^{A_n^*(s-r)}\big)A_{1}
^* y\big|_{H}ds\)^2\nonumber\\
&&\ds\q +  \(\int_{r}^T \big|e^{A^*_n(s-r)}
\big(A_{1}^*-A_{1,n}^*\big) y\big|_{H}ds\)^2+
 \(\int_r^T \big|e^{A^*_n(s-r)} A_{1,n}^*
(y-y_n)\big|_{H}ds\)^2\nonumber\\
&&\ds\q  + \(\int_r^T
\big|\big(e^{A^*(s-r)}-e^{A_n^*(s-r)}\big)C^*z\big|_{H}ds\)^2
+ \(\int_r^T \big|e^{A^*_n(s-r)}
\big(C^*-C_{n}^*\big)z\big|_{H}ds\)^2\\
&&\ds\q +  \(\int_r^T \big|e^{A^*_n(s-r)}
C_{n}^* (z-z_n)\big|_{H}ds\)^2 + \mE\(\int_r^T
\big|\big(e^{A^*(s-r)}-e^{A_n^*(s-r)}\big)Q^*x\big|_{H}ds\)^2
\nonumber\\
&&\ds\q + \(\int_r^T \big|e^{A^*_n(s-r)}
\big(Q^*-Q_{n}^*\big)x\big|_{H}ds\)^2+
\(\int_r^T \big|e^{A^*_n(s-r)}
Q_{n}^* (x-x_n)\big|_{H}ds\)^2\nonumber\\
&&\ds\q +  \int_r^T
\big|\big(e^{A^*(s-r)}-e^{A_n^*(s-r)}\big)z(s)\big|_{H}^2ds+
\int_r^T
\big|e^{A_n^*(s-r)}(z(s)-z_n(s))\big|_{H}^2ds\] \nonumber\\
&&\ds\q+4(14\wt M_T^2\!+\!1)\mE
\[\int_{T_1}^T\!
\big|\big(e^{A^*(T-s)}\!-\!e^{A_n^*(T-s)}\big)l(s)\big|_{H}^2
ds+ \!\int_{T_1}^T\!
\big|e^{A_n^*(T-s)}\big(l(s)\!-\!l_n(s)\big)\big|_{H}^2
ds\nonumber\\
&&\ds\q + \int_{T_1}^T\Big|\int_s^T
\big(e^{A^*(\si-s)}
\!-\!e^{A_n^*(\si-s)}\big)\k(\si,s)d\si\Big|_{H}^2ds\nonumber\\
&&\ds\q + \int_{T_1}^T\Big|\int_s^T
e^{A_n^*(\si-s)}
(\k(\si,s)-\k_n(\si,s))d\si\Big|_{H}^2ds\].\nonumber
\end{eqnarray}
Since
$$
\big|\big(e^{A^*(T-r)}-e^{A_n^*(T-r)}\big)G
x(T)\big|_{H}^2 \leq \cC\big|x(T)\big|_{H}^2,
$$
it follows from \eqref{6.8-eq5.2} and Lebesgue's
dominated convergence theorem that
\begin{equation}\label{2.9-eq36}
\begin{array}{ll}\ds
\lim_{n\to\infty}\mE\sup_{r\in
[T_1,T]}\big|\big(e^{A^*(T-r)}-e^{A_n^*(T-r)}\big)G
x(T)\big|_{H}^2=0.
\end{array}
\end{equation}
Similarly, we have that
\begin{equation}\label{2.9-eq37}
\lim_{n\to\infty}\mE\sup_{r\in
[T_1,T]}\big|e^{A_n^*(T-r)}(G-G_n)x(T)\big|_{H}^2=0.
\end{equation}
Similar to the proof of \eqref{2.9-eq18}, we can
obtain that
\begin{eqnarray}\label{2.9-eq38.1}
&&\ds \3n\3n\3n\3n\lim_{n\to\infty}\mE\sup_{r\in
[T_1,T]}\Big\{\[\int_r^T
\big|\big(e^{A^*(s-r)}-e^{A_n^*(s-r)}\big)A_{1}
^* y\big|_{H}ds\]^2 +  \[\int_r^T
\big|e^{A^*_n(s-r)} \big(A_{1}^*-A_{1,n}^*\big)
y\big|_{H}ds\]^2\nonumber\\
&&\ds\qq  +  \[\int_r^T
\big|\big(e^{A^*(s-r)}-e^{A_n^*(s-r)}\big)C^*z\big|_{H}ds\]^2
+  \[\int_r^T \big|e^{A^*_n(s-r)}
\big(C^*-C_{n}^*\big)z\big|_{H}ds\]^2 \nonumber\\
&&\ds\qq +  \[\int_r^T
\big|\big(e^{A^*(s-r)}-e^{A_n^*(s-r)}\big)Q^*x\big|_{H}ds\]^2
+  \[\int_r^T \big|e^{A^*_n(s-r)}
\big(Q^*-Q_{n}^*\big)x\big|_{H}ds\]^2\\
&&\ds \qq + \int_r^T
\big|\big(e^{A^*(s-r)}-e^{A_n^*(s-r)}\big)z\big|_{H}^2ds+
\int_{T_1}^T \big|\big(e^{A^*(T-s)}
-\!e^{A_n^*(T-s)}\big)l(s)\big|_{H}^2\! ds\nonumber\\
&&\qq \ds +  \int_{T_1}^T\!\Big| \int_s^T
\big(e^{A^*(\si-s)}
\!-\!e^{A_n^*(\si-s)}\big)\k(\si,s)d\si\Big|_{H}^2\!ds\Big\}\!\!=\!0.\nonumber
\end{eqnarray}

By the first equality in \eqref{6.7-eq3}, we get
that
\begin{equation}\label{2.9-eq38}
\begin{array}{ll}\ds
\lim_{n\to\infty}\mE\sup_{r\in
[T_1,T]}\[\big|e^{A_n^*(T-r)}G_n
(x(T)-x_n(T))\big|_{H}^2 + \(\int_r^T
\big|e^{A^*_n(s-r)} Q_{n}^*
(x-x_n)\big|_{H}ds\)^2\]\\
\ns\ds \leq \cC\lim_{n\to\infty}\[
\mE\big|x(T)-x_n(T)\big|_{H}^2 +
\(\int_{T_1}^T|Q^*|_{\cL(H)}ds\)^2\sup_{r\in
[T_1,T]}\mE \big|x(r)-x_n(r)\big|_{H}\]=0.
\end{array}
\end{equation}
Clearly,
\begin{equation}\label{2.9-eq39}
\begin{array}{ll}\ds
\mE\sup_{r\in [T_1,T]}\(\int_r^T
\big|e^{A^*_n(s-r)} A_{1,n}^*
(y-y_n)\big|_Hds\)^2+ \mE\sup_{r\in
[T_1,T]}\(\int_r^T
\big|e^{A^*_n(s-r)} C_{n}^* (z-z_n)\big|_Hds\)^2\\
\ns\ds \leq \!\wt M_T^2 \(\!\int_{T_1}^T\! \big|
A_{1}^*\big|_{\cL(H)}ds\)^2 \mE\!\!\sup_{r\in
[T_1,T]}\!|y(r)\!-\!y_n(r)|_H^2\!+\!\wt M_T^2
\(\int_{T_1}^T\!\! \big| C^*\big|_{\cL(H)}^2ds\)
\mE\!\int_{T_1}^T \!|z(r)\!-\!z_n(r)|_H^2ds.
\end{array}
\end{equation}

From \eqref{2.9-eq22} and \eqref{2.9-eq27}, we
have that
\begin{eqnarray}\label{2.9-eq32}
&&\ds\mE \int_{T_1}^T
\big|e^{A_n^*(T-s)}\big(l(s)-l_n(s)\big)\big|_H^2
ds\leq  \wt M_T^2\mE \int_{T_1}^T
\big|(l(s)-l_n(s))\big|_H^2 ds \nonumber\\
&&\ds \leq  \wt M_T^2\mE \big|G x(T)-\mE(G
x(T)|\;\cF_{T_1})-G_n x_n(T)+\mE(G_n
x_n(T)|\;\cF_{T_1})\big|_H^2\\
&&\ds \leq 2\wt M_T^2\mE \big|G x(T) -G_n x_n(T)
\big|_H^2 + 2\wt M_T^2\mE \big|\mE(G
x(T)|\;\cF_{T_1})-\mE(G_n
x_n(T)|\;\cF_{T_1})\big|_H^2\nonumber\\
&&\ds \leq 4\wt M_T^2\mE \big|G x(T) -G_n x_n(T)
\big|_H^2\leq 8\wt M_T^2\mE \big|(G-G_n) x(T)
\big|_H^2+\cC\wt M_T^2\mE \big|x(T) - x_n(T)
\big|_H^2.\nonumber
\end{eqnarray}
Since
$$
\big|(G-G_n) x(T) \big|_H^2\leq \cC\big|x(T)
\big|_H^2,
$$
it follows from Lebesgue's dominated convergence
theorem and \eqref{6.8-eq5} that
\begin{equation}\label{2.9-eq33}
\begin{array}{ll}\ds
\lim_{n\to\infty}\mE \big|(G-G_n) x(T)
\big|_H^2=\mE \lim_{n\to\infty}\big|(G-G_n) x(T)
\big|_H^2=0.
\end{array}
\end{equation}
By the first equality in \eqref{6.7-eq3} again,
we get that
$$
\lim_{n\to\infty}\mE \big|x(T) -  x_n(T)
\big|_H^2=0.
$$
This, together with \eqref{2.9-eq32} and
\eqref{2.9-eq33}, implies that
\begin{equation}\label{2.9-eq34}
\begin{array}{ll}\ds
\lim_{n\to\infty} \mE \int_{T_1}^T
\big|e^{A_n^*(T-s)}\big(l(s)-l_n(s)\big)\big|_H^2
ds =0.
\end{array}
\end{equation}
From \eqref{2.9-eq23}, \eqref{2.9-eq24},
\eqref{2.9-eq28} and \eqref{2.9-eq29}, we
conclude that
\begin{eqnarray}\label{2.9-eq35}
&&\3n\3n\ds\q 4(14\wt M_T^2+1)\mE
\int_{T_1}^T\Big|\int_s^T
e^{A_n^*(\si-s)}(\k(\si,s)-\k_n(\si,s))d\si\Big|_H^2ds\nonumber\\
&&\3n\3n\ds\leq 4(14\wt M_T^2+1)\wt
M_T^2\[\int_{T_1}^T \big(\mE|A_{1}^* y+ C^*z+
Qx-A_{1,n}^* y_n- C_n^*z_n-
Q_nx_n|_H^2\big)^{\frac{1}{2}}ds
\]^2\nonumber\\
&&\3n\3n\ds\leq\! 24(14\wt M_T^2\!+\!\!1)\wt
M_T^2
\Big\{\[\!\int_{T_1}^T\!\!\big(\mE|(A_1^*\!\!-\!A_{1,n}^*)y|_H^2\big)^{\frac{1}{2}}ds\]^2
\!\!+
\!\!\(\!\int_{T_1}^T\!\!|A_{1,n}^*|_{\cL(H)}
ds\)^2\mE\!\! \sup_{r\in
[T_1,T]}\!|y(r)\!-\!y_n(r)|_H^2\nonumber\\
&&\3n\3n\ds \q +
\[\int_{T_1}^T\big(\mE|(C^*-C_{n}^*)z|_H^2\big)^{\frac{1}{2}}ds\]^2
+ \(\int_{T_1}^T |C_{n}^*|_{\cL(H)}^2 ds\) \mE
\int_{T_1}^T |z-z_{n}|_H^2 ds\\
&&\3n\3n\ds \q +
\[\int_{T_1}^T\big(\mE|(Q^*-Q_{n}^*)x|_H^2\big)^{\frac{1}{2}}ds\]^2
+ \(\int_{T_1}^T |Q_{n}^*|_{\cL(H)} ds\)^2 \mE
\sup_{r\in[T_1,T]}|x(r)-x_{n}(r)|_H^2\Big\}.\nonumber
\end{eqnarray}
Combining \eqref{2.9-eq40},  \eqref{2.9-eq36},
\eqref{2.9-eq37}, \eqref{2.9-eq38},
\eqref{2.9-eq38.1}, \eqref{2.9-eq39},
\eqref{2.9-eq34} and \eqref{2.9-eq35}, we find
that
\begin{equation*}\label{2.9-eq30}
\begin{array}{ll}\ds
\lim_{n\to\infty}\(\mE\sup_{r\in
[T_1,T]}|y(r)-y_n(r)|_H^2 +
\mE\int_{T_1}^T|z(r)-z_n(r)|_H^2dr\)=0.
\end{array}
\end{equation*}
By repeating the above argument, we obtain the
second and third equality in \eqref{6.7-eq3}.
\endpf


\subsection{Proof of Propositions \ref{prop3}--\ref{prop2}}


{\it Proof of Proposition \ref{prop3}}\,: Let
$\tau\in [0,T)$ such that
$$
\int_\tau^T \big(2\big|A_1+B\Th\big|_{\cL(H)}
+\big|C+D\Th\big|_{\cL(H)}^2\big) ds
<\frac{1}{2M_T^2}.
$$
Define a map $\cG:C_\cS([\tau,T];\cL(H))\to
C_\cS([\tau,T];\cL(H))$ as follows:
$$
\begin{array}{ll}\ds
\cG(P)(r)\zeta\3n&\ds=e^{(T-r)A^*}Ge^{(T-r)A}\eta
+ \int_r^T
e^{(s-r)A^*}\big[P(A_1+B\Th)+(A_1+B\Th)^*
P\\
\ns&\ds\q+(C+D\Th)^* P(C+D\Th)+\,\Th^* R\Th
+Q\big]e^{(s-r)A}\zeta ds,\qq \forall\zeta\in H.
\end{array}
$$

Let $P_1,P_2\in C_\cS([0,T];\cL(H))$. Then, for
each $\zeta\in H$,
\begin{eqnarray*} &&\ds
\sup_{r\in[\tau,T]}\big|\big(\cG(P_1)-\cG(P_1)\big)(r)\zeta\big|_H\\
&&\ds= \sup_{r\in[\tau,T]}\Big|\int_r^T
e^{(s-r)A^*}[(P_1-P_2)(A_1+B\Th)+(A_1+B\Th)^*
(P_1-P_2)\\
&&\ds\qq\qq+(C+D\Th)^*
(P_1-P_2)(C+D\Th)]e^{(s-r)A}\zeta ds\Big|_H\\
&&\ds \leq \int_\tau^T
\big|e^{(s-r)A^*}\big|_{\cL(H)}\big[\big(\big|A_1+B\Th\big|_{\cL(H)}+\big|(A_1+B\Th)^*\big|_{\cL(H)}\big)\big|(P_1-P_2)\big|_{\cL(H)}
\\
&&\ds\qq\qq+\big|(C+D\Th)^*\big|_{\cL(H)}
\big|P_1-P_2\big|_{\cL(H)}\big|(C+D\Th)\big|_{\cL(H)}\big]\big|e^{(s-r)A}\big|_{\cL(H)}\big|\zeta\big|_H
ds\\
&&\ds\leq M_T^2\int_\tau^T
\big(2\big|A_1+B\Th\big|_{\cL(H)}
+\big|C+D\Th\big|_{\cL(H)}^2\big) ds
\sup_{r\in[\tau,T]}\big|(P_1-P_2)(r)\big|_{\cL(H)}
\big|\zeta\big|_H\\
&&\ds \leq
\frac{1}{2}\sup_{r\in[\tau,T]}\big|(P_1-P_2)(r)\big|_{\cL(H)}
\big|\zeta\big|_H.
\end{eqnarray*}
Therefore,
\begin{equation}\label{2.9-eq41}
\begin{array}{ll}\ds
\sup_{r\in[\tau,T]}\big|\big(\cG(P_1)-\cG(P_1)\big)(r)\big|_{\cL(H)}\leq
\frac{1}{2}\sup_{r\in[\tau,T]}\big|(P_1-P_2)(r)\big|_{\cL(H)}.
\end{array}
\end{equation}
This deduces that $\cG$ is contractive.
Consequently, there is a unique fixed point of
$\cG$, which is the mild solution to
\eqref{eq-Lya} on $[\tau,T]$. Repeating this
process gives us the unique $P\in
C_\cS([0,T];\cL(H))$ which satisfies
\eqref{8.20-eq22}. On the other hand, it is
obvious that if $P\in C_\cS([0,T];\cL(H))$
satisfying \eqref{8.20-eq22}, then $P^*\in
C_\cS([0,T];\cL(H))$ does. Hence, we have
$P=P^*$, that is, $P\in C_\cS([0,T];\cS(H))$.

The uniqueness of the solution is obvious.

\vspace{0.1cm} From \eqref{8.20-eq22}, we see
that for any $\zeta\in H$,
$$
\begin{array}{ll}\ds
\sup_{r\in[0,T]}\big|P(r)\zeta\big|_H\\
\ns\ds\leq
\sup_{r\in[0,T]}\big|e^{(T-r)A}\big|_{\cL(H)}^2\big|G\big|_{\cL(H)}
\big|\zeta\big|_H \!+\! \int_0^T\!
\big|e^{(s-r)A}\big|_{\cL(H)}^2
\big[\big(2\big|A_1+B\Th\big|_{\cL(H)}
\\
\ns\ds\q+\big|C+D\Th\big|_{\cL(H)}^2
\big)\big|P\big|_{\cL(H)}
+\big|\Th\big|_{\cL(H;U)}^2 \big|R\big|_{\cL(U)}
+\big|Q\big|_{\cL(H)}\big] \big|\zeta\big|_H ds.
\end{array}
$$
This, together with Gronwall's inequality,
implies the desired result.
\endpf

\vspace{0.3cm}

{\it Proof of Proposition \ref{prop1}} : For any
$\eta,\zeta\in H$, we have that
$$
\begin{array}{ll}\ds
\lan P(r)\eta,\zeta\ran\3n&\ds=\lan
Ge^{(T-r)A}\eta, e^{(T-r)A}\zeta\ran + \int_r^T
\lan [P(A_1+B\Th)+(A_1+B\Th)^*
P\\
\ns&\ds\q+(C+D\Th)^* P(C+D\Th)+\,\Th^* R\Th
+Q]e^{(s-r)A}\eta, e^{(s-r)A}\zeta\ran ds.
\end{array}
$$
If $\eta,\zeta\in D(A)$, it follows that $\lan
P(r)\eta,\zeta\ran$ is differentiable with
respect to $r$. A simple computation implies
\eqref{8.20-eq20}.
\endpf

\vspace{0.3cm}

The proof of Proposition \ref{prop2} is similar
to the one of Proposition \ref{prop1}, we omit
it.


\subsection{Proof Proposition \ref{prop4.1}}


The proof of Proposition \ref{prop4.1} is based
on a standard argument involving a minimizing
sequence and locally weak compactness of Hilbert
spaces. We give it in the appendix for the
completeness of the paper and the convenience
for some readers.

{\it Proof Proposition \ref{prop4.1}}\,: It
follows from the uniform convexity of
$u(\cd)\mapsto \cJ(0,0;u(\cd))$ that
\begin{equation}\label{J>l*.1}
\cJ(0,0;u(\cd))\geq\l
\dbE\1n\int_0^T|u(s)|^2ds,\qq\forall
u(\cd)\in\cU[0,T],
\end{equation}
for some $\l>0$. For any $t\in[0,T)$, and any
$u(\cd)\in\cU[t,T]$, let
\begin{equation}\label{ext}
v(s)\=\left\{\2n\ba{ll}0,\qq\ s\in[0,t),\\
\ns\ds u(s),\q
s\in[t,T].\ea\right.\end{equation}
Then $v(\cd) \in\cU[0,T]$. Since the initial
state is $0$, the solution $x(\cd)$ to
$$\left\{\2n\ba{ll}
\ns\ds dx =\big[(A+A_1)x +B v \big]ds+\big(C x +D v \big)dW(s) &\mbox{ in }(0,T], \\
\ns\ds x(0)=0,\ea\right.$$
satisfies
$$x(s)=0,\qq s\in[0,t].$$
Therefore,
\begin{equation}\label{4.4}
\cJ(t,0;u(\cd))=\cJ(0,0; v(\cd))\geq \l
\dbE\1n\int_0^T\big|v(s)\big|^2ds=\l
\dbE\1n\int_t^T|u(s)|^2ds.\end{equation}
Thus, $u(\cd)\mapsto \cJ(t,0;u(\cd))$ is
uniformly convex for any given $t\in[0,T)$. Let
$\bar u(\cd)\in\cU[t,T]$ be an optimal control
for $(t,\eta)$.

We can rewrite the cost functional as follows:
\begin{equation}\label{J-rep}
\begin{array}{ll}
\ns\ds \cJ(t,\eta;u(\cd))\\
\ns\ds=\lan G\big(\widehat\G_t\eta+\widehat
\Xi_tu\big),\widehat\G_t\eta+\widehat \Xi_tu\ran
+\lan Q(\G_t\eta+\Xi_tu),\G_t\eta+\Xi_tu\ran
+\langle Ru,u\rangle\\
\ns \ds=\lan\big(\widehat \Xi_t^*G\widehat
\Xi_t+\Xi_t^*Q\Xi_t+R\big)u,u\ran+2\lan\big(\widehat
\Xi_t^*G\widehat\G_t+\Xi_t^*Q\G_t\big)\eta,u\ran\\
\ns \ds\q
+\lan\big(\widehat\G_t^*G\widehat\G_t+\G_t^*Q\G_t\big)\eta,\eta\ran\\
\ns \ds = \cJ(t,\eta;0)+\cJ(t,0;
u(\cd))+2\int_t^T\lan\big(\widehat
\Xi_t^*G\widehat\G_t+\Xi_t^*Q\G_t\big)\eta,u\ran
ds\\
\ns\ds\geq  \l \dbE\int_t^T|\bar
u(s)|^2ds+\cJ(t,\eta;0)-\!\frac{\l}{2}\dbE\!\int_t^T\!\!|\bar
u(s)|^2ds\!
-\!\frac{1}{2\l}\dbE\!\int_t^T\!\!|\big(\widehat
\Xi_t^*G\widehat\G_t+\Xi_t^*Q\G_t\big)\eta|^2ds\\
\ns\ds\geq \frac{\l}{2}\dbE\int_t^T|\bar
u(s)|^2ds+\cJ(t,\eta;0)
-\frac{1}{2\l}\dbE\int_t^T| \big(\widehat
\Xi_t^*G\widehat\G_t+\Xi_t^*Q\G_t\big)\eta|^2ds.
\end{array}
\end{equation}
This implies that $u(\cd)\mapsto
\cJ(t,\eta;u(\cd))$ is coercivity. Clearly,
$u(\cd)\mapsto \cJ(t,\eta;u(\cd))$ continuous
and convex. Consequently, it has a unique
minimizer.

Moreover,  \eqref{J-rep} implies that
\begin{equation}\label{uni-convex-prop2}
V(t,\eta)\geq
\cJ(t,\eta;0)-\frac{1}{2\l}\dbE\int_t^T|\big(\widehat
\Xi_t^*G\widehat\G_t+\Xi_t^*Q\G_t\big)\eta|^2ds.
\end{equation}
Since the functions on the right-hand side of
(\ref{uni-convex-prop2}) are quadratic in $x$
and continuous in $t$, we get
(\ref{uni-convex-prop0}).
\endpf


\subsection{Proof of Proposition \ref{prop4.5}}


{\it Proof of Proposition \ref{prop4.5}} : Let
$\Th(\cd)\in L^2(0,T;\cL(U;H))$ and $P(\cd)$ be
the solution to \eqref{eq-Lya}. For any
$u(\cd)\in\cU[0,T]$, let $x_0(\cd)$ be the
solution to
$$
\left\{\2n
\begin{array}{ll}
\ns\ds dx_0 =\big[(A+A_1+B\Th) x_0+Bu\big]ds+\big[(C+D\Th)x_0+Du\big]dW(s) & \mbox{ in }(0,T], \\
\ns\ds x_0(0)=0.
\end{array}
\right.
$$
It follows from {\rm(\ref{J>l*})} and Lemma
\ref{lm2.3} that
$$
\l\dbE\!\int_0^T\!\!|\Th x_0 \!+\!u |^2ds
\!\leq\! \cJ(0,0;\Th(\cd)x_0(\cd)\!+\!u(\cd))
\!=\!\dbE\!\int_0^T\!\!\big[2\lan\big(L\!+\!K\Th\big)x_0,u\ran\!+\!\lan
Ku,u\ran\big]ds.
$$
Consequently, for any $u(\cd)\in\cU[0,T]$, it
holds that
\begin{equation}\label{P>LI}
\begin{array}{ll}
\ds\dbE\int_0^T\Big\{2\lan\big[L+(K-\l
I)\Th\big]x_0,u\ran+\lan (K-\l
I)u,u\ran\Big\}ds=\l\dbE\int_0^T|\Th(s)x_0(s)|^2ds\geq0.
\end{array}
\end{equation}
Fix any $u_0\in U$, take
$u(s)=u_0\chi_{[t,t+h]}(s)$, with $0\leq
t<t+h\leq T$. Then
\begin{equation}\label{8.20-eq35}
|x_0|_{C_\dbF([0,T];L^2(\Om;H))}\leq
\cC|u|_{L^2_\dbF(0,T;U)}\leq
\cC\sqrt{h}|u_0|_{U}.
\end{equation}
Dividing both sides of \eqref{P>LI} by $h$ and
letting $h\to 0$, noting \eqref{8.20-eq35}, we
obtain
$$\lan\big(K(s)-\l I\big)u_0,u_0\ran\geq 0,\qq\ae~s\in[0,T],\q \forall u_0\in U.$$
This gives the first inequality in
(\ref{Convex-prop-1}). To prove the second one,
for any $(t,\eta)\in[0,T)\times H$ and
$u(\cd)\in\cU[t,T]$, let $x_1(\cd)$ be the
solution to
$$
\left\{\2n
\begin{array}{ll}
\ds dx_1 =\big[(A+A_1+B\Th) x_1+Bu\big]ds+\big[(C+D\Th)x_1+Du\big]dW(s) &\mbox{ in }(t,T], \\
\ns\ds x_1(t)=\eta.
\end{array}
\right.
$$
It follows from Proposition \ref{prop4.1} and
Lemma \ref{lm2.3} that
$$
\begin{array}{ll}
\ds \a|x|^2\leq V(t,\eta)
\leq \cJ(t,\eta;\Th(\cd)x_1(\cd)+u(\cd))\\
\ns\ds =\langle
P(t)\eta,\eta\rangle\1n+\dbE\int_t^T\big[
2\lan\big(L+K\Th\big)x_1,u\ran\1n
+Ku,u\ran\big]ds.
\end{array}
$$
In particular, by taking $u(\cd)=0$ in the
above, we obtain
$$\lan P(t)\eta,\eta\ran \geq\a|x|^2,\qq\forall (t,\eta)\in[0,T]\times H,$$
which gives the second inequality.
\endpf


\subsection{Proof of Lemmas \ref{lm2.4} and \ref{lm2.5}}


{\it Proof of Lemma \ref{lm2.4}}:  Similar to
the proof of Proposition \ref{prop3}, one can
show the existence of a unique solution $P\in
C_\cS([0,T];\cS(H))$ to \eqref{lm2.4-eq1}. Let
$\tilde x(\cd)$ be the solution to the following
SEE:
$$
\left\{
\begin{array}{ll}\ds
d\tilde x =(A+\wt A )\tilde xdt+\wt C \tilde x
dW(t) &\mbox{ in }(t,T],\cr
\ns\ds \tilde x(t)=\eta\in H.
\end{array}
\right.
$$
Clearly this equation admits a unique solution.
By It\^o's formula  and Proposition \ref{prop1},
$$
\lan \wt P \eta,\eta\ran =\lan \wt G \tilde
x(T),\tilde x(T)\ran + \int_t^T \lan \wt
Q(s)\tilde x(s),\tilde x(s)\ran ds.
$$
This, together with \eqref{lm2.4-eq2}, implies
that $P(\cd)\in C([t,T];\cl{\dbS_+}(H))$.
\endpf

\ms

{\it Proof of Lemma \ref{lm2.5}}:  Let
$\Th(\cd)\in L^2(t,T;\cL(H;U))$. Define a
bounded linear operator
$\mathfrak{L}:\cU[t,T]\to\cU[t,T]$ by
$$\mathfrak{L}u=u-\Th x.$$
Then $\mathfrak{L}$ is bijective and its inverse
$\mathfrak{L}^{-1}$ is given by
$$\mathfrak{L}^{-1}u=u+\Th\wt x,$$
where $\wt x(\cd)$ is the solution to
$$
\left\{\2n
\begin{array}{ll}
\ns\ds d\tilde x =\big[\big(A+A_1 +B \Th \big)\wt x +B u \big]ds+ \big[\big(C +D \Th \big)\wt x +D u \big]dW(s) &\mbox{ in }(t,T], \\
\ns\ds\wt x(t)=0.
\end{array}
\right.
$$
By the bounded inverse theorem,
$\mathfrak{L}^{-1}$ is bounded with norm
$|\mathfrak{L}^{-1}|_{\cL(\cU[t,T])}>0$. Thus,
$$\begin{array}{ll}
\ns\ds\dbE\int_t^T|u(s)|_U^2ds
=\dbE\int_t^T|(\mathfrak{L}^{-1}\mathfrak{L}u)(s)|_U^2ds
\leq|\mathfrak{L}^{-1}|\dbE\int_t^T|(\mathfrak{L}u)(s)|_U^2ds\\
\ns\ds\qq\qq\qq~~\1n=|\mathfrak{L}^{-1}|_{\cL(\cU[t,T])}\dbE\int_t^T\big|u(s)-\Th(s)
x(s)\big|_U^2ds, \qq\forall
u(\cd)\in\cU[t,T],\end{array}$$
which implies (\ref{lem-2.6}) with
$c_0=|\mathfrak{L}^{-1}|_{\cL(\cU[t,T])}^{-1}$.
\endpf

\end{document}